\documentclass[11pt]{article}

\usepackage[margin=1in]{geometry}
\usepackage{mathtools,amsfonts,amssymb,mathrsfs,graphicx,bm}
\usepackage{verbatim}
\usepackage{stmaryrd}
\usepackage{algorithm,algpseudocode}
\usepackage{siunitx}
\usepackage{booktabs}
\usepackage[hidelinks]{hyperref}
\usepackage[capitalize,nameinlink]{cleveref}

\title{A Mean-Informed Low-Rank Monolithic Stochastic Galerkin Solver for the Unsteady Navier--Stokes Equations}

\author{Ahmet Kaan Aydin\thanks{Department of Mathematics and Statistics, University of Maryland, Baltimore County, 1000 Hilltop Circle, Baltimore, MD 21250, USA.
  (\texttt{aaydin1@umbc.edu}, \texttt{sousedik@umbc.edu}).}
\and Bed\v{r}ich Soused\'{\i}k\footnotemark[1]}
\date{}

\DeclareMathOperator{\diag}{diag}

\hypersetup{
  pdftitle={A Mean-Informed Low-Rank Monolithic Stochastic Galerkin Solver for the Unsteady Navier--Stokes Equations},
  pdfauthor={A. Aydin, and B. Soused\'{\i}k}
}

\begin{document}

\maketitle

\begin{abstract}

We study a low-rank solver for the stochastic unsteady incompressible Navier--Stokes equations with uncertain viscosity. The problem is discretized by a stochastic Galerkin method and written in an all-at-once (monolithic) form where the time steps are informed by a sequential solve of the mean problem. The solution vector is represented in Tensor Train (TT) format, while the system matrices are represented in CANDECOMP/PARAFAC (CP) format as needed. We propose preconditioners based on the low-rank CP approximations that retain additional stochastic and nonlinear information compared to mean-based preconditioners. Numerical experiments for a benchmark channel-flow problem are presented to show the effectiveness of the low-rank approximation, the relative tolerance strategy, and the proposed preconditioners.
\end{abstract}

\noindent\textbf{Keywords:} uncertainty quantification, spectral stochastic finite element method, Navier--Stokes equations, stochastic Galerkin method, low-rank representation, tensor decomposition.

\noindent\textbf{MSC 2020:} 35R60, 60H15, 65N22, 65N30, 65N35.

\section{Introduction}
\label{sec:intro}
\noindent Physical systems with uncertainty are often modeled by stochastic partial differential equations. 
In this paper, a stochastic version of the incompressible Navier--Stokes equations in which the viscosity is modeled as a random field is studied. 
Consequently, velocity and pressure become random fields too.
Uncertainty in viscosity may arise from measurement errors, contaminants, multiphase mixtures, or external effects. 
Such models also appear in applications including blood flow~\cite{pereira2013uncertainty}, aeroelasticity~\cite{pettit2004uncertainty}, and nuclear engineering~\cite{puscas2023interaction}. 
We use a spectral stochastic Galerkin approach based 
on generalized polynomial chaos (gPC) expansions~\cite{ghanem2003stochastic,le2010spectral,xiu2010numerical,xiu2002wiener}, 
and seek the corresponding gPC expansions of the velocity and pressure. 

Many traditional solvers for the deterministic problem solve the system sequentially in time~\cite{elman2011fast, elman2014finite}, 
and only some of them allow for adaptive selection of the step size~
\cite{kay2010adaptive,sousedik2016stochastic,sousedik2022stochastic}. 
Here, we also first use a sequential solver for the underlying deterministic problem with the mean viscosity 
to predetermine the time-step sizes~\cite{sousedik2022stochastic}. 
Subsequently, in the formulation of the stochastic problem, the equations at these time steps are assembled into an all-at-once formulation. 
The coupling through both the time-steps and the terms of gPC expansion leads to a single large system. 
This poses on one hand a challenge. 
However, looking at the entire stochastic Galerkin problem at once offers an opportunity to leverage its somewhat hidden structure 
and address the complexity of the problem by means of low-rank tensor methods~\cite{bachmayr2023low,kolda2009tensor,grasedyck2013literature}. 
Specifically, the solution vector is represented in Tensor Train (TT) format~\cite{oseledets2011tensor}. 
Such strategy 
was developed in~\cite{elman2020low}, but it turns out that the nonlinear convection term can still cause rapid growth of TT ranks. 
To address this issue, we represent the convection and stochastic operators using CANDECOMP/PARAFAC (CP) format with low-rank approximations 
computed by Alternating Least Squares (ALS)~\cite{carroll1970analysis,harshman1970foundations}. 
We note that a further motivation for combining the all-at-once formulation with low-rank representations arises from the growing gap between processor speed and memory performance~\cite{mutlu2022modern, patterson1997case}. The proposed approach reduces communication costs via the all-at-once formulation and compresses storage through low-rank tensor decompositions. Together, these steps help shift the computational bottleneck away from memory access, enabling more efficient large-scale simulations.

Our contributions thus advances methods for the stochastic unsteady Navier--Stokes equations~
\cite{elman2020low}. In addition to extending the use of low-rank tensor methods and nonuniform time-stepping, 
we also implement a flexible tolerance strategy for TT-rounding and GMRES iterations. 
The CP approximations of the convection and stochastic operators improve the efficiency of the solver, 
and reduce the rank growth caused by the matrix-vector products in TT format. 
The CP approximations offer the same low-rank operators as the mean-based preconditioner, but retain more 
information. 
This motivated the proposal of a CP-based preconditioner. 
Numerical experiments 
illustrate the effectiveness of the low-rank approximation, the relative tolerance strategy, and the proposed preconditioners.

The rest of the paper is organized as follows. \Cref{sec:unsteady-ns} presents the deterministic unsteady incompressible Navier--Stokes equations. \Cref{sec:stochastic-galerkin} shows stochastic Galerkin formulation, and the resulting all-at-once discrete system. \Cref{sec:solution} describes the nonlinear 
iteration and the preconditioned Krylov method used to solve the linearized systems. \Cref{sec:low-rank-approximation} introduces the TT representation of the solution, the CP-ALS approximation of the operators, the preconditioner construction, and the flexible tolerance strategy. \Cref{sec:numerical_SGTT} gives numerical results 
and examines the performance of the monolithic solver. \Cref{sec:conclusion_1} summarizes the main findings.

\section{Unsteady Deterministic Navier--Stokes Equations}
\label{sec:unsteady-ns}
\noindent Consider the unsteady incompressible Navier--Stokes~\cite{elman2014finite, temam2024navier} problem in a spatial domain $D\subset \mathbb{R}^2$ given by
\[
\begin{array}{rl}
    \displaystyle\frac{\partial \vec u}{\partial t} - \nu \nabla^2 \vec u + (\vec u \cdot \nabla) \vec u + \nabla p = \vec f &\\
    \nabla \cdot \vec u = 0  &\text{ in } D \times (0, T],
\end{array}
\]
where $\vec u$ is the velocity field, $p$ is the pressure, and $\nu$ is the kinematic viscosity. 
The boundary conditions are Dirichlet boundaries~$\Gamma_{\rm Dir}$ for an inflow and no-slip on the walls, 
and Neumann boundary set for natural outflow on $\Gamma_{\rm Neu}$ given, respectively, by
\[
\vec u = \vec g  \text{ on } \Gamma_{\rm Dir} \times (0, T], \qquad     \nu \nabla \vec u \cdot \vec n - p \vec n = \vec 0   \text{ on } \Gamma_{\rm Neu} \times (0, T].
\]
Consider the temporal discretization by implicit (backward) Euler method 
\begin{align}
\frac{\vec{u}^{\,k}-\vec{u}^{\,k-1}}{\tau_{k}}-\nabla\cdot\left(  \nu\nabla\vec
{u}^{\,k}\right)  +\vec{u}^{\,k}\cdot\nabla\vec{u}^{\,k}+\nabla p^{\,k}  &  =\vec
{f}^{\,k},\label{eq:NS2-a}\\
\nabla\cdot\vec{u}^{\,k}  &  =0, \label{eq:NS2-b}%
\end{align}
\noindent where $\tau_k$ denotes the time step size, while $\vec{u}^{\,k}$ and $p^{\,k}$ represent the velocity and pressure at time $t_k$, respectively. The spatial discretization is achieved using div-stable finite elements. 
In numerical experiments we use Taylor--Hood finite element discretization and quadrilateral elements. 
Define the finite-dimensional spaces
\begin{align}
	\mathcal{X}=\operatorname*{span}\{\bm{\phi	}_{i}(x)\}_{i=1}^{n_u}\subset\left(  H^{1}(D)\right)  ^{2},\qquad
	\mathcal{Y}=\operatorname*{span}\{\varphi_{i}(x)\}_{i=1}^{n_p}\subset L^{2}(D).
\end{align}
Let $\mathcal{X}_{\text{Dir}}^{k}$ be the space of functions in $\mathcal{X}$ with Dirichlet boundary conditions~$\vec{u}_{\text{Dir}}$. The finite
element solutions, $\vec{u}^{\,k}\in\mathcal{X}_{\text{Dir}}^{k}$ and $p^{\,k}\in\mathcal{Y}$, are approximated as
\begin{equation}\label{eq:discrete-solutions}
	\vec{u}^{\,k}=\sum_{i=1}^{n_{u}}{u}_{i}^{\,k}\bm{\phi}_{i}(x),\qquad {p}^{\,k}%
=\sum_{i=1}^{n_{p}}{p}_{i}^{\,k}{\varphi}_{i}(x),\qquad k=1,\dots,n_{t}.
\end{equation}
The coefficient vectors $\bm{u}^{k}=[{u}_{1}^{k},\dots,{u}_{n_{u}}^{k}]$ and
$\bm{p}^{k}=[{p}_{1}^{k},\dots,{p}_{n_{p}}^{k}]$ are computed from
\begin{align}
\tau_{k}^{-1} \bm M \bm{u}^{k}-\tau_{k}^{-1}\bm M \bm{u}^{k-1}+\bm A \bm{u}^{k}+\bm N(\bm{u}^{k})\bm{u}^{k}%
+\bm B^{T}\bm{p}^{k}  &  =\bm{f}^{u,k},\label{eq:NS3-a}\\
\bm B \bm{u}^{k}  &  =\bm{f}^{p,k}, \label{eq:NS3-b}%
\end{align}
where $\bm M $\ is the velocity mass matrix, $\bm A$ is the vector-Laplacian
matrix, $\bm N(\bm{u}^{k})$ is the vector-convection matrix, and $\bm B$ is the divergence matrix. These matrices are defined as
\begin{align*}
	\bm M   &  =\left[  m_{cd}\right]  ,\quad m_{cd}=\int_{D}\bm{\phi}_{d}\,\bm{\phi}_{c}, &
\bm A  &  =\left[  a_{cd}\right]  ,\quad a_{cd}=\nu\int_{D}\nabla\bm{\phi}_{d}%
:\nabla\bm{\phi}_{c},\\
\bm N(\bm{u}^{k})  &  =\left[  n_{cd}\right]  ,\quad n_{cd}=\int_{D}(\bm{u}^{k}\cdot
\nabla\bm{\phi}_{d})\bm{\phi}_{c},&\bm B  &  =\left[  b_{cd}\right]  ,\quad b_{cd}=-\int_{D}{\varphi}_{c}\left(
\nabla\cdot\bm{\phi}_{d}\right).
\end{align*}
The boundary conditions are incorporated into the right-hand side of~(\ref{eq:NS3-a})--(\ref{eq:NS3-b}), which may be written as
\begin{equation}
\left[
\begin{array}
[c]{cc}%
\bm{F}_{u}^{k} & \bm{B}^{T}\\
\bm{B} & 0
\end{array}
\right]  \left[
\begin{array}
[c]{c}%
\bm{u}^{k}\\
\bm{p}^{k}%
\end{array}
\right]  +\left[
\begin{array}
[c]{cc}%
-\tau_{k}^{-1}\bm{M} & 0\\
0 & 0
\end{array}
\right]  \left[
\begin{array}
[c]{c}%
\bm{u}^{k-1}\\
\bm{p}^{k-1}%
\end{array}
\right]  =\left[
\begin{array}
[c]{c}%
\bm{f}^{u,k}\\
\bm{f}^{p,k}%
\end{array}
\right], \label{eq:NS4}
\end{equation}
where $\bm{F}_{u}^{k}=\tau_{k}^{-1}\bm{M}+\bm{A}+\bm{N}(\bm{u}^{k}),$ for $k=1,\ldots, n_{t}.$

\section{The Stochastic Galerkin Formulation} \label{sec:stochastic-galerkin}

\noindent Let $\left(  \Omega,\mathcal{F}_{\Omega},\mathcal{P}\right)  $ be a complete probability space, where$~\Omega$ is the sample space, $\mathcal{F}_{\Omega}$ is the $\sigma$-algebra of measurable events on$~\Omega$, and $\mathcal{P}$ is a probability
measure. Assume that the stochasticity in the viscosity is induced by a
vector of independent random variables $\bm \xi=(\xi_{1},\dots,\xi_{m_{\xi}})^{T}$
such that $\xi:\Omega\rightarrow\Gamma\subset%
\mathbb{R}
^{m_{\xi}}$. Let$~\mathcal{B}(\Gamma)$ denote the Borel $\sigma$-algebra
on$~\Gamma$ induced by$~\xi$, and let$~\mu$ denote the induced measure. The
expected value of the product of two measurable functions of$~\xi$ on$~\Gamma$
determines a Hilbert space $T_{\Gamma}\equiv L^{2}(\Omega,F_{\xi},\mu)$ with
inner product
\[
\left\langle u,v\right\rangle =\mathbb{E}\left[  uv\right]  =\int_{\Gamma
}u(\xi)v(\xi)\,d\mu(\xi),
\]
where the symbol $\mathbb{E}$ denotes mathematical expectation. This infinite-dimensional space is approximated by a finite-dimensional polynomial subspace $\mathcal{Z}\subset
T_{\Gamma}$ spanned by a set of polynomials $\left\{  \psi_{\ell}%
(\xi)\right\}  $ that are orthogonal with respect to the density
function $\mu$, that is $\mathbb{E}\left[  \psi_{k}\psi_{\ell}\right]
=\delta_{k\ell}$. These functions form the generalized polynomial chaos (gPC) basis, see~\cite{ghanem2003stochastic,xiu2002wiener} for details. The dimension of the space $\mathcal{Z}$ depends on the polynomial degree. For polynomials of total
degree $p$, the dimension is $n_{\xi}=\binom{m_{\xi}+p}{p}$.
Specifically, the viscosity is given by a gPC expansion as
\begin{equation}
\label{eq:gPC-nu}
\nu(x,\xi)=\sum_{\ell=1}^{n_{\nu}}\nu_{\ell}(x)\,\psi_{\ell}(\xi),
\end{equation}
where$~\nu_{1}$ is the mean viscosity ($\psi_{1}\equiv1$) and $\nu_{\ell}$,
$\ell=2,\dots,n_{\nu}$, are a set of given deterministic functions, 
obtained for example from a Karhunen-Lo\`{e}ve expansion. 
At each time step, the velocity and pressure~\eqref{eq:discrete-solutions} are expanded in the gPC basis as
\begin{equation}
\vec{u}^{\,k}=\sum_{\ell=1}^{n_{\xi}}\sum_{i=1}^{n_{u}}u_{i\ell}^{k}\bm{\phi
}_{i}(x)\psi_{\ell}(\xi),~~~ p^{k}=\sum_{\ell=1}^{n_{\xi}}\sum_{i=1}%
^{n_{p}}p_{i\ell}^{k} \varphi_{i}(x)\psi_{\ell}(\xi),~~~ k=1,\dots,n_{t}.
\label{eq:gPC-expansions}%
\end{equation}
Substituting~\eqref{eq:gPC-nu}--\eqref{eq:gPC-expansions} into~\eqref{eq:NS4} and applying stochastic Galerkin projection yields 
\begin{equation*}
		\left[
\begin{array}
[c]{cc}%
\bm{F}_{u,\rm{SG}}^{k} & \bm{I}_{n_\xi}\otimes{\bm B}^{T}\\
\bm{I}_{n_\xi}\otimes\bm{B} & \bm 0
\end{array}
\right]  \left[
\begin{array}
[c]{c}%
\bm{u}^{k}\\
\bm{p}^{k}%
\end{array}
\right] -\left[
\begin{array}
[c]{cc}%
\bm{I}_{n_\xi}\otimes\tau_{k}^{-1}\bm{M} & \bm 0\\
\bm 0 & \bm 0
\end{array}
\right]  \left[
\begin{array}
[c]{c}%
\bm{u}^{k-1}\\
\bm{p}^{k-1}%
\end{array}
\right]  =\left[
\begin{array}
[c]{c}%
\bm{f}^{u,k}\\
\bm{f}^{p,k}%
\end{array}
\right],%
\end{equation*}
where $\bm{I}_{n_\xi}$ denotes the identity matrix of size $n_\xi$ and $\otimes$ denotes the Kronecker product. 
The matrix $\bm{F}_{u,\rm{SG}}^{k}$ is defined as
\begin{equation}
\bm{F}_{u,\rm{SG}}^{k}=\bm{I}_{n_{\xi}}\otimes\tau_{k}^{-1}\bm{M}+\sum_{\ell=1}^{n_{\nu}}\left(  \bm H_{\ell}\otimes
\bm A_{\ell}\right) +\sum_{\ell=1}^{n_{\xi}}\bm H_{\ell
}\otimes \bm N(u_{\ell}^{k}), \qquad k=1,\ldots, n_{t}, \label{eq:F_u^k_stoch}%
\end{equation}
where
\[
    \bm A_{\ell}    =\left[  a_{\ell,cd}\right]  ,~ a_{\ell,cd}=\int_{D}\nu
_{\ell}(x)\,\nabla\phi_{d}:\nabla\phi_{c}, \quad
\bm N(u_{\ell}^{k})    =\left[  n_{\ell,cd}\right]  ,~ n_{\ell,cd}=\int%
_{D}(u_{\ell}^{k}\cdot\nabla\phi_{d})\phi_{c},
\]
and the stochastic matrices $\bm H_{\ell}$ are given by
\[
\bm H_{\ell}=\left[ \bm  H_{\ell}\right]  _{ij},\quad\left[ \bm H_{\ell}\right]
_{ij}=\mathbb E\left[  \psi_{\ell}\psi_{i}\psi_{j}\right]  ,\quad i,j=1,\dots n_{\xi
},\quad\ell=1,\dots, n_\nu.
\]

\subsection{All-at-once Formulation}
\label{subsec:all_at_once}
We consider the time steps concatenated as
\[
\mathbf{u}=\left[
\begin{array}
[c]{c}%
\mathbf{u}^{1}\\
\vdots\\
\mathbf{u}^{n_{t}}%
\end{array}
\right]  ,\quad\mathbf{u}^{k}=\left[
\begin{array}
[c]{c}%
\bm u_{1}^{k}\\
\vdots\\
\bm u_{n_{\xi}}^{k}%
\end{array}
\right]  ,\quad \bm u_{\ell}^{k}=\left[
\begin{array}
[c]{c}%
u_{1\ell}^{k}\\
\vdots\\
u_{n_{u}\ell}^{k}%
\end{array}
\right]  \text{,}\quad\text{etc.}, 
\]
which corresponds to the operators defined as 
\begin{align}
\mathbb{B} &:= \bm I_{n_t}\otimes \bm I_{n_\xi}\otimes \bm B, &
\mathbb{M} &:= \bm T\otimes \bm I_{n_\xi}\otimes \bm M,\label{eq:matrices1}\\[4pt]
\mathbb{A} &:= \bm I_{n_t}\otimes\!\left(\sum_{\ell=1}^{n_\nu} \bm H_\ell\otimes \bm A_\ell\right), &
\mathbb{N}(\bm u) &:= \mathop{\mathrm{blkdiag}}\limits_{k=1,\dots,n_t}\left(\sum_{\ell=1}^{n_\xi} \bm H_\ell\otimes \bm N(\mathbf{u}^k_\ell)\right),\label{eq:matrices2}
\end{align}
where, with~$\bm D=\diag(\tau_1^{-1},\tau_2^{-1},\dots,\tau_{n_t}^{-1})$, the matrix $\bm T$ is
\begin{equation}\label{eq:TimeMatrixDeterministic}
\bm T
=
\bm D\bm E,
~
\bm E=
\begin{bmatrix}
1 &  &  & \\
-1 & 1 &  & \\
& \ddots & \ddots & \\
&  & -1 & 1
\end{bmatrix}\text{ and } \bm{E}^{-1} = \begin{bmatrix}
1 & 0 & \dots & 0 \\
1 & 1 & \dots & 0 \\
\vdots & \vdots & \ddots & \vdots \\
1 & 1 & \dots & 1
\end{bmatrix}.
\end{equation}
So, the convection block $\mathbb{N}(\bm{u})$ is block-diagonal in time (nonlinearity acts only within each time level) and couples stochastic modes through $H_\ell$. Note that the coupling among the time steps by backward Euler method is contained in $\mathbb M$ via the matrix~$\bm T$. With these definitions, the global saddle-point system reads
\begin{equation}
\label{eq:all-at-once}
\left[
\begin{array}
[c]{cc}%
\mathbb{F}_u & \mathbb{B}^{T}\\
\mathbb{B} & 0
\end{array}
\right]  \left[
\begin{array}
[c]{c}%
\mathbf{u}\\
\mathbf{p}%
\end{array}
\right]  =\left[
\begin{array}
[c]{c}%
\bm f^{u}\\
	\bm f^{p}%
\end{array}
\right],
\qquad
\mathbb{F}_u :=\mathbb{M}
+ \mathbb{A}
+ \mathbb{N}(\bm{u}),
\end{equation}
and the vectors $\bm f^u$ and $\bm f^p$ incorporate the initial and boundary conditions.

\section{Nonlinear Iteration}
\label{sec:solution}

\noindent We use Picard method for the nonlinear problem~\eqref{eq:all-at-once}.  
Each Picard step entails solving a linear saddle-point system, which is done inexactly using (right-)preconditioned GMRES. 
The preconditioner is based on approximate Schur complement. Its application entails solves with the block $F_u$ performed inexactly using (a nested) GMRES. 
We continue to use the notational convention: Bold lowercase letters denote vectors, bold uppercase letters denote matrices, blackboard bold letters denote all-at-once matrices and calligraphic symbols denote tensor objects. 
When the nonlinear convection term corresponds to Picard step~$i$, we write
\[
\mathbb{F}_u^{(i)} := \mathbb{F}_u(\bm{u}^{(i)}).
\]

\subsection{Picard's Method}
\label{subsec:picard_schur}
\noindent Let $\mathbf{u}^{(i)},\ \mathbf{p}^{(i)}$ be the approximate solution at the $i$th step of the iteration. Each Picard step entails solving the corresponding residual equation for a correction of the solution. Let
$\mathbf{u}^{(i)} = \mathbf{u}^{(i-1)} + \delta \mathbf{u}^{(i)}$ and $\mathbf{p}^{(i)} = \mathbf{p}^{(i-1)} + \delta \mathbf{p}^{(i)}$.
Then $\delta \mathbf{u}^{(i)}$ and $\delta \mathbf{p}^{(i)}$ satisfy
\begin{equation}
\left[
\begin{array}{cc}
\mathbb{F}_u^{(i-1)} & \mathbb{B}^{T} \\
\mathbb{B} & 0
\end{array}
\right]
\left[
\begin{array}{c}
\delta \mathbf{u}^{(i)} \\
\delta \mathbf{p}^{(i)}
\end{array}
\right]
=
\left[
\begin{array}{c}
\bm r^{u,(i-1)} \\
\bm r^{p,(i-1)}
\end{array}
\right],
\label{eq:Picard}
\end{equation}
where the nonlinear residual is
\[
\bm r^{(i)} =
\left[
\begin{array}{c}
\bm r^{u,(i)} \\
\bm r^{p,(i)}
\end{array}
\right]
=
\left[
\begin{array}{c}
\bm{f}^{u} \\
\bm{f}^{p}
\end{array}
\right]
-
\left[
\begin{array}{cc}
\mathbb{F}_u^{(i)} & \mathbb{B}^{T} \\
\mathbb{B} & 0
\end{array}
\right]
\left[
\begin{array}{c}
\mathbf{u}^{(i)} \\
\mathbf{p}^{(i)}
\end{array}
\right].
\]
The solution is updated 
until the prescribed tolerance of the residual $\bm r^{(i)}$ is met.  

\subsubsection{Preconditioner for a Picard step} 
\label{subsubsec:outer_gmres}

\noindent At each Picard iteration, residual equation corresponding to~\eqref{eq:Picard} is solved by a right-preconditioned flexible GMRES~\cite{saad2003iterative}.
The block upper-triangular preconditioner considered in this study is given by
\begin{equation}
\mathbb{P} = \left[
\begin{array}
[c]{cc}%
\mathbb{F}_{u}^{-1} & \mathbb{F}_{u}^{-1}\mathbb{B}^{T}\mathbb{S}^{-1}\\
0 & -\mathbb{S}^{-1}%
\end{array}
\right], \label{eq:prec}%
\end{equation}
where $\mathbb{S} \approx \mathbb{B}\mathbb{F}_{u}^{-1}\mathbb{B}^{T}$, i.e., it approximates the exact Schur complement. 
Note that, with a suitable ordering of the block operations, an application of~$\mathbb{P}$ requires only one application of $\mathbb{F}_{u}^{-1}$. 
Since use of the exact Schur complement is computationally prohibitive, 
its practical approximations are typically derived from the fact that the commutator of the convection-diffusion operator 
with the divergence operator is small under certain assumptions~\cite[Chapter 8]{elman2014finite}. 
Two standard methods that exploit this property are the pressure convection-diffusion (PCD) and the least-squares commutator (LSC) preconditioners.

A detailed comparison of the PCD and LSC methods in the context of low-rank stochastic Navier--Stokes solvers is provided in~\cite{elman2020low}. 
Here, we use the LSC preconditioner. 
The primary advantage of the LSC preconditioner is that it avoids the explicit construction of boundary conditions and convection-diffusion matrices on the pressure space. Instead, it operates directly on the existing velocity space matrices. The LSC approximation for the inverse Schur complement is given by 
\begin{equation}
\mathbb{S}^{-1}_{\text{LSC}} = (\mathbb{B}\widehat{\mathbb{M}}^{-1}\mathbb{B}^{T})^{-1}(\mathbb{B}\widehat{\mathbb{M}}^{-1}\mathbb{F}_{u}\widehat{\mathbb{M}}^{-1}\mathbb{B}^{T})(\mathbb{B}\widehat{\mathbb{M}}^{-1}\mathbb{B}^{T})^{-1},
\label{eq:schur-approx-lsc}
\end{equation}
where the action of the inverse of the velocity mass matrix, $\mathbb{M}^{-1}$, is approximated by $\widehat{\mathbb{M}}^{-1} := \bm T^{-1}\otimes \bm I_{n_\xi}\otimes \bm M_{*}^{-1}$, where subscript $*$ denotes the diagonal of the corresponding matrix. The inverse temporal matrix $\bm{T}^{-1} = \bm{E}^{-1}\bm{D}^{-1}$ is known explicitly by \eqref{eq:TimeMatrixDeterministic}.

\subsubsection{Application of \texorpdfstring{$\mathbb{F}_u^{-1}$}{F\_u\^(-1)}}
\label{subsubsec:inner_gmres}
\begin{figure}[ht]
\centering
\includegraphics[width=\linewidth]{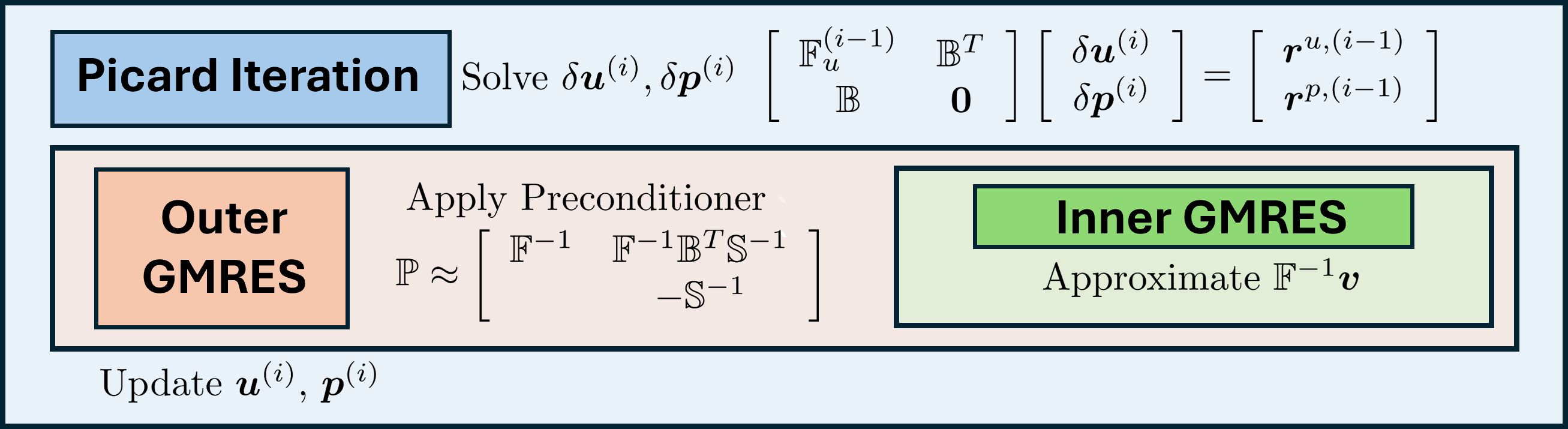}
    \caption{Schematic overview of the nested iterative solver. Each Picard step entails one solve using an outer GMRES method 
    and one solve as an action of $\mathbb{F}^{-1}$, which is approximated by an inner low-rank preconditioned GMRES solver.}
\label{fig:Diagram}
\end{figure}

\noindent The application of the outer preconditioner $\mathbb{P}$ in \eqref{eq:prec} requires the action of $\mathbb{F}_u^{-1}$, which is the computationally dominant part of the solver. For a given vector $\bm v$, the product $\mathbb{F}_u^{-1}\bm v$ is computed approximately by solving
$
\mathbb{F}_u \bm x = \bm v
$
with GMRES. This leads to a nested iterative strategy. At each Picard step, the residual equation is solved by a preconditioned GMRES method, referred to as the \emph{outer} GMRES. Each application of the outer preconditioner requires an approximate solve with $\mathbb{F}_u$, which is performed by a second GMRES method, referred to as the \emph{inner} GMRES. \Cref{fig:Diagram} provides a schematic overview of this strategy. 

Previous work by McDonald et al.~\cite{mcdonald2016simple} explored a block-diagonal preconditioner for this system, but it exhibited slow convergence, often requiring a number of iterations proportional to the number of time steps. Elman and Su~\cite{elman2020low} proposed a more effective preconditioner based on the mean solution, specifically designed for problems with a 
constant time step size. 
By utilizing the dimensionless temporal difference matrix $\bm{E}$ from~\eqref{eq:TimeMatrixDeterministic}, 
the mean-based preconditioner \cite[Eq. (5.23)]{elman2020low} is 
\begin{equation}\label{eq:innerGMRES_Elman}
\mathbb{P}_{\text{avg}} = \bm{E}^{-1} \otimes \bm{I}_{n_\xi} \otimes (\tau^{-1}\bm{M} + \bm{A}_{1} + \bm{N}(\bm{u}_{\text{avg}}))^{-1},
\end{equation}
where $\bm{A}_1$ is the mean vector-Laplacian matrix and the convection matrix $\bm{N}(\bm{u}_{\text{avg}})$ is evaluated at the mean solution. Unlike \eqref{eq:innerGMRES_Elman}, our formulation uses adaptive, nonuniform time stepping ($\tau_k$). To obtain a mean-based baseline for the numerical experiments in \cref{sec:numerical_SGTT}, we approximate the uniform step with the average, $\tau^{-1} \approx \bar{\tau}^{-1} = \frac{1}{n_t}\sum_{k=1}^{n_t}\tau_k^{-1}$. While this provides a workable baseline for the nonuniform regime, the approach exhibits reduced efficiency when there are large variations in time step sizes. The effectiveness of this strategy is discussed in Section~\ref{sec:numerical_SGTT}. The preconditioner~\eqref{eq:innerGMRES_Elman} disregards the stochastic coupling intrinsic to the system. Furthermore, the time coupling is artificially enforced on the vector-Laplacian and convection matrices. To address these limitations, we propose a novel preconditioner based on the low-rank representation of the $\mathbb{F}_u$ matrix. This structure-preserving preconditioner is described in detail in Subsection~\ref{subsubsec:Fu_inv}.

\section{Low-Rank Approximation}
\label{sec:low-rank-approximation}
\noindent The all-at-once system~\eqref{eq:all-at-once} exhibits a natural separability:
temporal, stochastic, and spatial operators appear (except for convection) as explicit
Kronecker products, see~\eqref{eq:matrices1}--\eqref{eq:matrices2}.
To exploit this structure, a TT representation for the 
velocity (and, when needed, pressure) vectors is considered. 
This choice enables matrix-vector products and preconditioner applications to be carried out by acting on each factor separately.

\subsection{Tensor Train Decomposition of the Solution Vector}
\label{subsec:tt_decom}
\noindent A three-way tensor $\mathcal{U}\in\mathbb{R}^{n_t\times n_\xi\times n_u}$ admits a \emph{TT-decomposition}, cf.~\cite{oseledets2011tensor}, as 
\[
\mathcal{U}(i_1,i_2,i_3)
\;\approx\;
\sum_{\alpha=1}^{\kappa_1}\sum_{\beta=1}^{\kappa_2}
u^{(1)}_{i_1,\alpha}\;
u^{(2)}_{\alpha,i_2,\beta}\;
u^{(3)}_{\beta,i_3}\,,
\]
where $\kappa_1,\kappa_2$ are the TT-ranks. The TT-decomposition is a low-rank representation of the tensor, which allows for efficient storage by reducing the memory requirements from $O(n_t n_\xi n_u)$ to $O(n_t \kappa_1 + n_\xi \kappa_1 \kappa_2 + n_u \kappa_2)$.
The tensor has an equivalent form in terms of Kronecker products, cf.~\cite{dolgov2014alternating}, 
\begin{equation}\label{eq:sol_vecform}
\text{vec}(\mathcal{U})
=\sum_{\alpha,\beta}
\bigl(\bm u^{(1)}_\alpha\bigr)\otimes
\bigl(\bm u^{(2)}_{\alpha,\beta}\bigr)\otimes
\bigl(\bm u^{(3)}_\beta\bigr).
\end{equation}
Any matrix $\mathbb X = \bm{X}^{(1)}\otimes\bm{X}^{(2)}\otimes\bm{X}^{(3)}$ with Kronecker product structure, such as matrices in \eqref{eq:matrices1}--\eqref{eq:matrices2}, its product with a vector of the form \eqref{eq:sol_vecform} can be computed efficiently as
\begin{equation}\label{eq:cheal_matvec}
\mathbb{X}\text{vec}(\mathcal{U}) = \sum_{\alpha,\beta}
\bigl(\bm{X}^{(1)} \bm u^{(1)}_\alpha\bigr)\otimes
\bigl(\bm{X}^{(2)} \bm u^{(2)}_{\alpha,\beta}\bigr)\otimes
\bigl(\bm{X}^{(3)} \bm u^{(3)}_\beta\bigr).
\end{equation}
The efficiency relies on maintaining low TT ranks, that is, keeping $\kappa_1$ and $\kappa_2$ modest. The TT ranks are preserved under matrix-vector products given that $\mathbb{X}$ has a rank-one Kronecker product structure. However, the sum of two TT tensors $\bm{u}_1$ and $\bm{u}_2$ with ranks $\kappa_1^{(1)}, \kappa_2^{(1)}$ and $\kappa_1^{(2)}, \kappa_2^{(2)}$ results in a tensor with increased rank,
\(
\kappa_1 = \kappa_1^{(1)} + \kappa_1^{(2)},~\kappa_2 = \kappa_2^{(1)} + \kappa_2^{(2)}.
\)
To control this growth, the TT representations truncated using the TT-SVD algorithm~\cite[Algorithm 1]{oseledets2011tensor} to a prescribed tolerance $\epsilon_{\text{TT}}$ such as
 \begin{equation}
\label{eq:relative-error}
\frac{\|\hat{\bm{u}} - \bm{u}\|_F}{\|\bm{u}\|_F} < \epsilon_{\text{TT}}.
\end{equation}
\noindent A summarized version of the algorithm for three-way tensors in~\Cref{alg:TT_SVD}.

\begin{algorithm}[ht]
\caption{\cite[Algorithm 1]{oseledets2011tensor} TT-SVD for Three-Way Tensor}
\label{alg:TT_SVD}
\begin{algorithmic}[1]
\Require Tensor $\mathcal{U}\in\mathbb{R}^{n_t\times n_\xi\times n_u}$, truncation tolerances $\epsilon_{\text{TT}}$
\Ensure TT cores $\{\bm{U}^{(1)},\bm{U}^{(2)},\bm{U}^{(3)}\}$

\State Set $\kappa_0 = 1$ and $\bm U \leftarrow \mathcal{U}$
\For{$j = 1,2$}
    \State Reshape $\bm U$ into matrix form
    $
    \bm U \in
    \begin{cases}
        \mathbb{R}^{n_t \times (n_\xi n_u)}, & j = 1,\\[2pt]
        \mathbb{R}^{(\kappa_{1} n_\xi) \times n_u}, & j = 2,
    \end{cases}
    $
    \State Compute truncated SVD $\bm U =\bm  W \bm{\Sigma} \bm V^{T}$ with $\|\bm U -\bm  W\Sigma \bm V^{T}\|_F \le \epsilon_{\text{TT}}$
    \State Set $\kappa_j = \operatorname{rank}(\Sigma)$
    \State Define TT core $\bm{U}^{(j)} = \text{reshape}(\bm W, [\kappa_{j-1}, n_j, \kappa_j])$
    \State Update $\bm U \leftarrow \bm{\Sigma} \bm V^{T}$
\EndFor
\State Set final core $\bm{U}^{(3)} = \bm U$
\State \Return TT representation of $\mathcal{U}$ with cores $\{\bm{U}^{(1)},\bm{U}^{(2)},\bm{U}^{(3)}\}$
\end{algorithmic}
\end{algorithm}

\subsection{CANDECOMP/PARAFAC (CP) Representation}
\label{subsec:cp-decomposition}
\noindent The rank $R$ CP representation of a three-way tensor $\mathbf{\mathcal{X}}\in\mathbb{R}^{n_t^2\times n_\xi^2\times n_u^2}$ is
\begin{equation*}
	\mathbf{\mathcal{X}} \approx \sum_{r=1}^R \lambda_r\, \bm{a}_r \circ \bm{b}_r \circ \bm{c}_r = \llbracket\bm{\lambda};\bm{A},\bm{B},\bm{C}\rrbracket,
\end{equation*}
where $\bm{A}\in\mathbb{R}^{n_t^2\times R}$, $\bm{B}\in\mathbb{R}^{n_\xi^2\times R}$, and $\bm{C}\in\mathbb{R}^{n_u^2\times R}$ are factor matrices, $\bm{\lambda}\in\mathbb{R}^R$ contains weights, and $\circ$ denotes the outer product. The CP decomposition is typically computed using the \emph{Alternating Least Squares} (ALS) algorithm shown in \Cref{alg:ALS_CP}, refer to \cite[Algorithm 1]{TTB_CPRALS} for details. We map the system operators (such as $\mathbb{N}$ and $\mathbb{A}$)  into an equivalent third-order CP tensor as
\begin{equation}
    \mathbb{X} = \sum_{m=1}^M \bm{X}_m^{(1)} \otimes \bm{X}_m^{(2)} \otimes \bm{X}_m^{(3)} \quad \to \quad \mathbf{\mathcal{X}} = \sum_{m=1}^M \text{vec}(\bm{X}_m^{(1)}) \circ \text{vec}(\bm{X}_m^{(2)}) \circ \text{vec}(\bm{X}_m^{(3)}).
    \label{eq:tensor_reshape}
\end{equation}
This reshaping yields a CP tensor of rank $M$. Then, the ALS algorithm (\Cref{alg:ALS_CP}) is applied to find a lower-rank approximation $\hat{\mathbf{\mathcal{X}}}$ with rank $R \ll M$.

\begin{algorithm}[ht]
\caption{\cite[Algorithm 1]{TTB_CPRALS} Alternating Least Squares (ALS)}
\label{alg:ALS_CP}
\begin{algorithmic}[1]
\Require Third-order tensor $\mathcal{X} \in \mathbb{R}^{n_1 \times n_2 \times n_3}$,
target CP rank $R$, tolerance $\varepsilon_{\mathrm{tol}}$, maximum iterations $k_{\max}$
\Ensure Factor matrices $\bm A \in \mathbb{R}^{n_1 \times R}$,
$\bm B \in \mathbb{R}^{n_2 \times R}$,
$\bm C \in \mathbb{R}^{n_3 \times R}$, and weights $\bm\lambda \in \mathbb{R}^R$

\State Initialize factor matrices $\bm A_0,\bm B_0,\bm C_0$ and weights $\bm\lambda$
\While{$r_k > \varepsilon_{\mathrm{tol}}$ \textbf{and} $k < k_{\max}$}
    \State $k \gets k+1$
    \State Form diagonal matrix $\bm\Lambda = \operatorname{diag}(\bm\lambda)$
    \State Update $\bm A_k, ~\bm B_k, ~\bm C_k$ by solving the least-squares problems
    \begin{align*}
			\bm{A}_k &= \arg\min_{\bm{A}}\bigl\|\bm{X}_{(1)} - \bm{A}\,\bm{\Lambda}\,(\bm{C}_{k-1}\odot \bm{B}_{k-1})^{T}\bigr\|_F^2,\\
			\bm{B}_k &= \arg\min_{\bm{B}}\bigl\|\bm{X}_{(2)} - \bm{B}\,\bm{\Lambda}\,(\bm{C}_{k-1}\odot \bm{A}_{k})^{T}\bigr\|_F^2,\\
			\bm{C}_k &= \arg\min_{\bm{C}}\bigl\|\bm{X}_{(3)} - \bm{C}\,\bm{\Lambda}\,(\bm{B}_{k}\odot \bm{A}_{k})^{T}\bigr\|_F^2,
		\end{align*}
        where $\odot$ denotes the Khatri-Rao product.
    \State Construct the rank-$R$ CP approximation
    $
    \widehat{\mathcal X}_k
    =
    \sum\limits_{r=1}^{R}\lambda_r\,\bm a_r^{(k)} \circ \bm b_r^{(k)} \circ \bm c_r^{(k)}
    $
    \State Compute the relative residual
    $
    r_k
    =
    \frac{\|\mathcal X-\widehat{\mathcal X}_k\|_F}{\|\mathcal X\|_F}
    $
\EndWhile
\State \Return $\bm A_k,\bm B_k,\bm C_k,\bm\lambda$
\end{algorithmic}
\end{algorithm}

An advantage of this reshaping strategy is the exact preservation of the Frobenius norm ($\|\mathbb{X}\|_F = \|\mathbf{\mathcal{X}}\|_F$). Thus, the residual condition is retained as
\begin{equation}
  r_k =   \frac{\|\mathbf{\mathcal{X}} - \hat{\mathbf{\mathcal{X}}}\|_F}{\|\mathbf{\mathcal{X}}\|_F} \equiv \frac{\|\mathbb{X} - \hat{\mathbb{X}}\|_F}{\|\mathbb{X}\|_F}.
\end{equation}
Once ALS converges to the desired tolerance in the tensor format, the optimized vector columns $\bm{a}_r, ~\bm{b}_r,$ and $\bm{c}_r$ are reshaped back into their respective matrix dimensions to assemble $\hat{\mathbb{X}}$. This structured reduction procedure is summarized in \Cref{alg:Operator_CP_ALS}.

\begin{algorithm}[ht]
\caption{Operator Rank Reduction via Tensorization and CP-ALS}
\label{alg:Operator_CP_ALS}
\begin{algorithmic}[1]
\Require Matrix $\mathbb{X} = \sum_{m=1}^M \bm{X}_m^{(1)} \otimes \bm{X}_m^{(2)} \otimes \bm{X}_m^{(3)}$, target rank $R$, tolerance $\varepsilon_{\text{tol}}$
\Ensure Low-rank matrix $\hat{\mathbb{X}} = \sum_{r=1}^R \hat{\bm{X}}_r^{(1)} \otimes \hat{\bm{X}}_r^{(2)} \otimes \hat{\bm{X}}_r^{(3)}$

\State \textbf{Tensorization:} Vectorize components to form initial CP tensor:
    \[ \mathbf{\mathcal{X}} = \sum_{m=1}^M \text{vec}(\bm{X}_m^{(1)}) \circ \text{vec}(\bm{X}_m^{(2)}) \circ \text{vec}(\bm{X}_m^{(3)}) \]
\State \textbf{CP-ALS:} Apply \cref{alg:ALS_CP} to tensor $\mathbf{\mathcal{X}}$ with rank $R$ and tolerance $\varepsilon_{\text{tol}}$.
\State \textbf{Reconstruction:} Reshape each vector back into its original matrix dimensions:
    \begin{equation*}
        \hat{\bm{X}}_r^{(i)} \leftarrow \text{reshape}(\bm{a}_r, \text{size}(\bm{X}_1^{(i)})), \quad i = 1,2,3
    \end{equation*}
\State \Return $\hat{\mathbb{X}} = \sum_{r=1}^R \hat{\bm{X}}_r^{(1)} \otimes \hat{\bm{X}}_r^{(2)} \otimes \hat{\bm{X}}_r^{(3)}$
\end{algorithmic}
\end{algorithm}
\noindent\subsection{Low-Rank Representations of Solutions and Convection Matrix}
\label{subsec:low-rank-representation}
The solution vector $\bm{u}$ can be represented in TT format as
\begin{equation}
	\label{eq:solution_ttform}
	\bm{u}
\;=\;
\sum_{\alpha=1}^{\kappa_1}\sum_{\beta=1}^{\kappa_2}\bigl(u^{(1)}_\alpha\bigr)\otimes
\bigl(u^{(2)}_{\alpha,\beta}\bigr)\otimes
\bigl(u^{(3)}_\beta\bigr).
\end{equation}
Similarly the pressure vector $\bm{p}$ will be considered to be represented in the TT format. This representation allows for efficient matrix-vector products as described in \eqref{eq:cheal_matvec}. Using the TT representation of the solution the convection matrix known to admit following Kronecker product form~\cite[Eq. 4.12]{elman2020low}
	\begin{equation}
	\label{eq:conv_matrix}
	\mathbb{N}(\bm u) =\sum_{\alpha,\beta}\left(  \text{diag}( u^{(1)}) \otimes \sum_{\ell=1}^{n_\xi} (u^{(2)}(\alpha,\ell,\beta)\bm H_\ell)  \otimes \bm N( u^{(3)}_{\beta})\right)
	\end{equation}
Note that the number of terms scales with the TT-ranks, \(\kappa_1\kappa_2\) of the solution and the stochastic dimension, \(n_\xi\). For simplicity, assume $\kappa_1 \approx \kappa_2 \approx \kappa$. The convection matrix $\mathbb{N}(\bm{u})$ in \eqref{eq:conv_matrix} consists of summation of $\mathcal{O}(\kappa^2 n_\xi)$ matrices, each written as a Kronecker product of three matrices. Recall that summation of two TT-tensors results in a TT-tensor with ranks equal to the sum of the individual ranks. Thus, multiplying $\mathbb{N}(\bm{u})$ by $\bm u$, in general any TT-formatted vector with comparable ranks to the solution $\bm u$, yields a tensor with ranks growing from $\mathcal{O}(\kappa)$ to $\mathcal{O}(n_\xi \kappa^3)$. Similarly, product of $\mathbb{A}$ with a TT-formatted vector with ranks $\kappa$ yields a tensor with $\mathcal{O}( n_\xi \kappa^2)$ terms, leading to TT rank growth of the form $\mathcal{O}(\kappa) \to \mathcal{O}(n_\xi\kappa)$.

The TT rounding (TT-SVD) is a crucial step in managing the rank growth of the solution state during the iterative solution process. However, the quadratic growth due to convection operator remains a challenge. Several strategies have been proposed to mitigate this rank growth by reducing the rank of the convection operator. These include replacing the convection matrix with one based on a low-rank approximation of the solution vector, or using the mean solution in its place~\cite{elman2020low}. While these approximations help maintain manageable TT ranks, they may fail to capture the full complexity of the system, potentially degrading the performance of the preconditioner or the accuracy of the solution.

The common theme in these approaches is to construct a lower-rank surrogate for $\mathbb{N}(\bm{u})$, that can be used in the preconditioners too. We consider the approximations of $\mathbb{N}(\bm{u})$ and $\mathbb{A}$ by low-rank CP decompositions. The approximation obtained by CP-ALS reduces the number of Kronecker summands from $\mathcal{O}(\kappa^2 n_\xi)$ (convection) and $\mathcal{O}(n_\xi)$ (Laplacian) to arbitrarily small CP rank, $R$, prescribed in the CP-ALS algorithm. Let the rank-$R$ CP approximations of the convection and diffusion matrices be
\begin{equation}
		\mathbb{N}(\bm{u}) \approx \hat{\mathbb{N}}(\bm{u}) = \sum_{r=1}^{R_1} \bm{N}_r^{(1)} \otimes  \bm{N}_r^{(2)} \otimes \bm{N}_r^{(3)},\quad
\mathbb{A}\approx \hat{\mathbb{A}} = \sum_{r=1}^{R_2} \bm{A}_r^{(1)} \otimes  \bm{A}_r^{(2)} \otimes \bm{A}_r^{(3)}.
\label{eq:cp-decomposition}
\end{equation}
The matrix–vector product involving $\hat{\mathbb{N}}(\bm{u})$ then results in a representation with $\mathcal{O}(R \kappa^2)$ terms, which is significantly smaller than the $\mathcal{O}(\kappa^4 n_\xi)$ terms required by the full convection matrix $\mathbb{N}(\bm{u})$. Moreover, rather than discarding particular dynamics from the system, this approach provides a controlled approximation of the operator itself, which can be tuned via the CP rank \(R\). This allows for a more balanced trade-off between computational efficiency and solution accuracy.

\subsubsection{Approximating \texorpdfstring{$\mathbb{F}_{u}^{-1}$}{F\_u\^{-1}}}
\label{subsubsec:Fu_inv}
The preconditioner~\eqref{eq:prec} requires applying $\mathbb{F}_{u}^{-1}$, which is computationally expensive. We therefore seek a structured low-rank surrogate for $\mathbb{F}_{u}^{-1}$ that preserves the Kronecker structure and is cheap to apply. One strategy is to construct a rank-one CP approximation of $\mathbb{F}_u$ directly.
\begin{equation}
    \mathbb{F}_{u} \approx \hat{\mathbb{F}}_{u} = \bm{F}_1^{(1)} \otimes \bm{F}_1^{(2)} \otimes \bm{F}_1^{(3)}.
\end{equation}
Because this surrogate consists of a single Kronecker product, its exact inverse is computationally trivial to apply as
\begin{equation}
    \mathbb{P}_{\text{CP1}}^{-1} = \left(\bm{F}_1^{(1)}\right)^{-1} \otimes \left(\bm{F}_1^{(2)}\right)^{-1} \otimes \left(\bm{F}_1^{(3)}\right)^{-1}.
    \label{eq:prec-cp1-inv}
\end{equation}
This approach successfully captures the dominant directional dynamics of the coupled flow while ensuring that the necessary matrix-vector products during the inner GMRES iterations remain highly efficient and tightly bounded in TT-rank.

Although \eqref{eq:prec-cp1-inv} makes the preconditioner cheap to apply, a rank-one approximation may not be sufficiently accurate for all problems. We propose a generalization to arbitrary-rank CP approximations
\begin{align}
	\hat{\mathbb{F}}_{u} =& \bm{T} \otimes \bm{I}_{n_\xi} \otimes \bm{M} + \hat{\mathbb{A}} + \hat{\mathbb{N}}(\bm{u})\nonumber\\ \label{eq:low-rank-F_u}
	\approx& (\bm{T} \otimes \bm{I}_{n_\xi} \otimes \bm{M} ) \left(\bm{I}+ \sum_{r=1}^{R_1} \bm{T}^{-1}\bm{A}_r^{(1)} \otimes  \bm{A}_r^{(2)} \otimes  \bm{M}^{-1}_*\bm{A}_r^{(3)}\right.\\ \nonumber
	&\hspace{3cm}\left. + \sum_{r=1}^{R_2}  \bm{T}^{-1}\bm{N}_r^{(1)} \otimes  \bm{N}_r^{(2)} \otimes  \bm{M}_{*}^{-1}\bm{N}_r^{(3)} \right)
\end{align}
The matrices $\bm{T}^{-1}$ and $\bm{M}_{*}^{-1}$ are known explicitly and can be applied efficiently. For small time steps, the term $\bm{T}\otimes \bm{I}_{n_\xi}\otimes \bm{M}$ becomes dominant in $\hat{\mathbb{F}}_{u}$. The factorization in~\eqref{eq:low-rank-F_u} isolates this dominant block and collects the remaining couplings into a perturbation. For brevity, denote
\begin{equation}
\mathbb{K}\;=\; \sum_{r=1}^{R_1} \bm{T}^{-1}\bm{A}_r^{(1)} \otimes  \bm{A}_r^{(2)} \otimes  \bm{M}_{*}^{-1}\bm{A}_r^{(3)}
\;+\; \sum_{r=1}^{R_2}  \bm{T}^{-1}\bm{N}_r^{(1)} \otimes  \bm{N}_r^{(2)} \otimes  \bm{M}_{*}^{-1}\bm{N}_r^{(3)}.
\end{equation}
Substituting this into \eqref{eq:low-rank-F_u} and inverting yields
\begin{equation}
\hat{\mathbb{F}}_{u}^{-1}= (\mathbb{I}+ \mathbb{K})^{-1}\,(\bm{T} \otimes \bm{I}_{n_\xi} \otimes \bm{M} )^{-1}.
\end{equation}

A natural approach is to approximate $(\mathbb{I}+\mathbb{K})^{-1}$ by a truncated Neumann series. However, convergence is not guaranteed in general, since the presence of $\bm{T}^{-1}$ in $\mathbb{K}$ may help reduce the spectral radius but does not ensure that $\mathbb{K}$ is small in norm. Therefore, an artificial damping factor $\delta\in(0,1]$ is introduced to ensure $
\rho(\delta\mathbb{K})<1.
$
With this damping, $(\mathbb{I}+\delta\mathbb{K})^{-1}$ can be approximated by a truncated Neumann series, giving the preconditioner
\begin{equation}\label{eq:prec-neumann-expansion}
    \mathbb{P}_{\mathrm{Neu}}^{-1}
    =
    \sum_{n=0}^{K}(-\delta \mathbb{K})^n
    \,(\bm{T}\otimes \bm{I}_{n_\xi}\otimes \bm{M}_*)^{-1},
\end{equation}
where $K$ is the truncation order. A heuristic choice of $\delta$ may be sufficient, but it can fail when the perturbation is not small. In such cases, the size of $\mathbb{K}$ can be estimated more directly using
\[
\rho(\mathbb{K})
\le
\tau_{\max}(1+r_k)
\bigl(\|\mathbb{A}\|_F+\|\mathbb{N}\|_F\bigr)
\|\bm{M}_*^{-1}\|_2.
\]
Here, $r_k$ measures the relative error in the CP-ALS, and $\tau_{\max}$ is the largest time step.
\subsection{Relative Rounding Tolerances and Overall Low-Rank Algorithm}
\label{subsec:relative-tols}

The overall structure of the low-rank algorithm is provided in Algorithm~\ref{alg:picard_krylov}. Each Picard iteration, \eqref{eq:Picard}, is solved by a LSC preconditioned GMRES. Each application of $\mathbb F_u^{-1}$ is computed by inner GMRES with a low-rank preconditioner, \eqref{eq:prec-cp1-inv} or \eqref{eq:prec-neumann-expansion}. Setting a priori tolerances for the rounding of the TT-decompositions and for the inner/outer GMRES solves can cause oversolving in early Picard steps. To avoid this, we propose a strategy that ties each iterations tolerance to the corresponding residual of the overarching iteration. All Krylov vectors and intermediate products are stored in TT format and truncated after each \texttt{matvec} by \emph{relative} rounding that depends on the relative residuals of the corresponding iterate. Thus, TT ranks remain modest in early Picard/outer stages and increase as required.
\begin{algorithm}[ht]
\caption{All-at-once Low-rank Picard-Krylov Solver}
\label{alg:picard_krylov}
\begin{algorithmic}[1]
\Require Initial guess $(\bm u^{(0)},\bm p^{(0)})$, tolerances $\varepsilon_{\rm Pic}$, $\varepsilon_{\rm out}$, $\varepsilon_{\rm in}$, $\varepsilon_{\rm TT}$, and $i_{\max}$
\Ensure Approximate solution $(\bm u,\bm p)$

\State Compute $\bm r^{(0)}$

\While{$\|\bm r^{(i)}\|/\|\bm r^{(0)}\| > \varepsilon_{\rm Pic}$ \textbf{and} $i < i_{\max}$}
    \State Form $\mathbb F_u^{(i)} = \mathbb F_u(\bm u^{(i)})$
    \State Compute CP surrogate $\widehat{\mathbb N}(\bm u^{(i)})$
    \State Solve the Picard correction $(\delta\bm u^{(i)},\delta\bm p^{(i)})$ with tolerance
    $
    \eta_{\rm out}^{(i)} = \varepsilon_{\rm out}\,\frac{\|\bm r^{(i)}\|}{\|\bm r^{(0)}\|}
    $
    \State In each application of $\mathbb P$, approximate $\mathbb F_u^{-1}$ with tolerance
    $
    \eta_{\rm in}^{(i,j)} = \varepsilon_{\rm in}\,\frac{\|\bm r_{\rm out}^{\,j}\|}{\|\bm r_{\rm out}^{\,0}\|}
    $
    and TT-rounding tolerance proportional to the current inner residual
    \State Update
    $
    (\bm u^{(i+1)},\bm p^{(i+1)})
    =
    (\bm u^{(i)},\bm p^{(i)})+(\delta\bm u^{(i)},\delta\bm p^{(i)})
    $
    \State Compute $\bm r^{(i+1)}$, $i \gets i+1$
\EndWhile

\State \Return $(\bm u,\bm p)=(\bm u^{(i)},\bm p^{(i)})$
\end{algorithmic}
\end{algorithm}

\section{Numerical Experiments}
\label{sec:numerical_SGTT}

\begin{figure}[ht]
    \centering
    \includegraphics[width=0.45\linewidth]{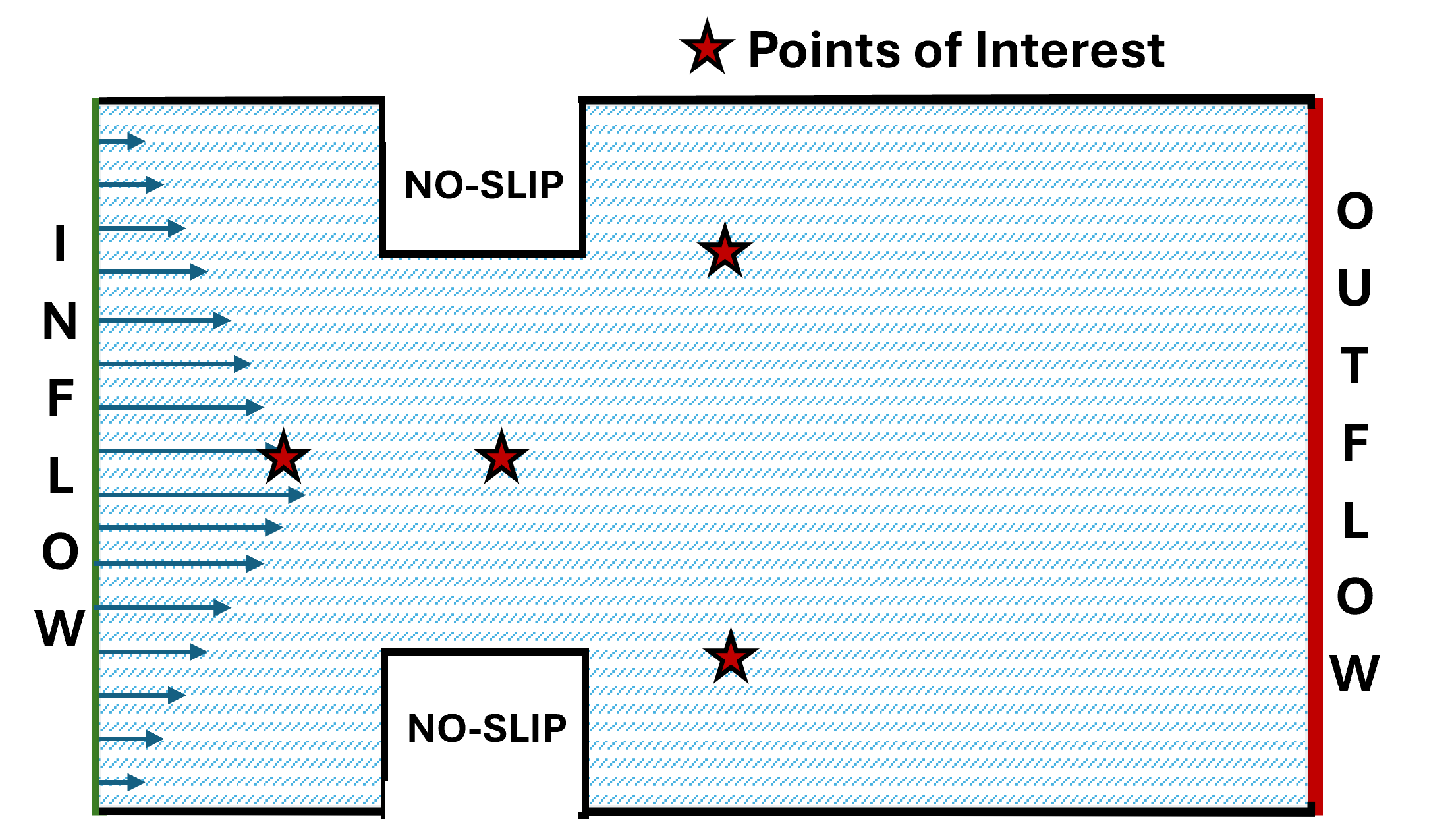} \includegraphics[width=0.4\linewidth]{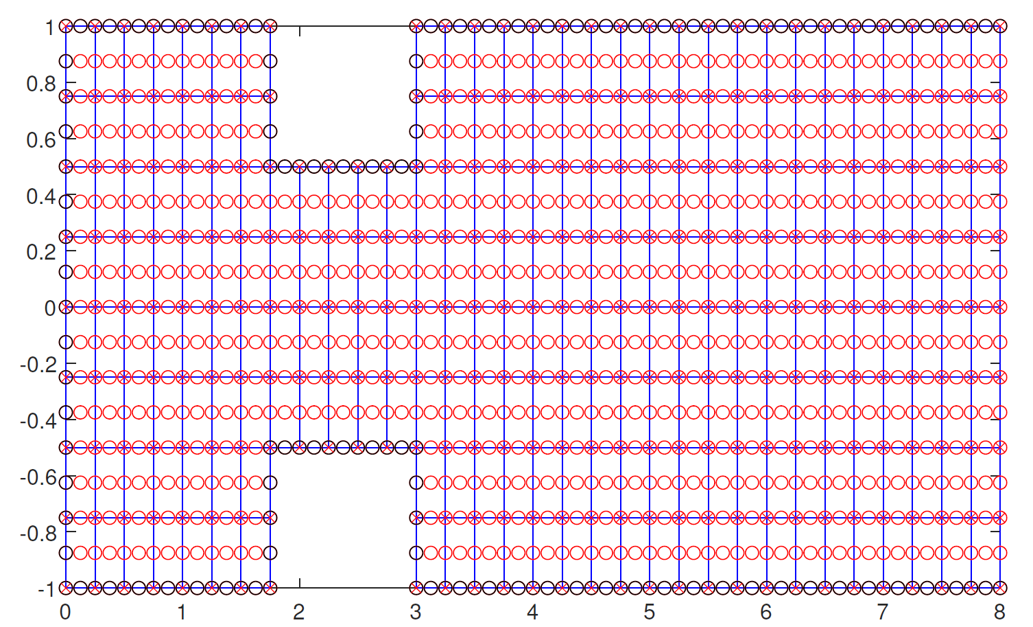}
    \caption{Narrow Channel flow domain with points used for probability density estimates (left) and the corresponding quadrilateral mesh (right). A parabolic inflow is prescribed at the left boundary, no-slip conditions on the walls, and natural outflow at the right boundary. Stars mark pressure degrees of freedom ($Q_1$), and circles mark velocity degrees of freedom ($Q_2$).}
    \label{fig:narrow_channel}
\end{figure}

Consider the flow through a Narrow Channel, a well-known benchmark problem widely used to assess the robustness of saddle-point preconditioners~\cite{elman2014finite, elman2020low}. The physical domain $D$ consists of a rectangular channel of length $L = 8$ and height $H = 2$, featuring a symmetrically constricted middle section where the height reduces to $1$. The boundary conditions are prescribed as follows: a time-dependent parabolic inflow profile is imposed at the left boundary ($x=0$), no-slip boundary conditions ($\bm{u} = \bm{0}$) are enforced on the top and bottom walls, and a natural outflow condition is applied at the right boundary. The parabolic inflow is
\begin{equation}\label{eq:inflow}
\bm{u}_{\mathrm{in}}(y, t) = \begin{bmatrix}
    1-y^2\\ 
    0
\end{bmatrix}\left(1-e^{-10t}\right).
\end{equation}
\noindent The spatial discretization is performed using Taylor-Hood finite elements on a uniform quadrilateral mesh, generated via the MATLAB IFISS package, version 3.7 \cite{elman2014ifiss, silvester2016ifiss}. The generated mesh is illustrated in \cref{fig:narrow_channel}. Numerical experiments were performed in MATLAB R2025b on a 64-bit Linux with Intel Xeon E5-2620 and 1.0 TiB of RAM. TT and CP operations are performed using the TT-Toolbox \cite{oseledets2011tensor, tt-toolbox} and the Tensor Toolbox \cite{tensor-toolbox}, respectively.

\subsection{Stochastic Viscosity Modelling}
\label{subsec:stoch_visc}

The kinematic viscosity $\nu$ is modeled as a truncated lognormal process with a specified mean value $\nu_1 = 0.01$ \cite{ghanem1999nonlinear, sousedik2016stochastic}. Let $G(x;\bm\xi)$ be a zero-mean Gaussian random field parameterized by a set of $N$ independent standard normal random variables $\bm\xi=(\xi_1,\dots,\xi_N)$. The strictly positive viscosity field is defined pointwise by
\begin{equation}\label{eq:nu-lognormal-ours}
    \nu(x;\bm\xi)=\exp\big(\mu_\nu + G(x;\bm\xi)\big) \quad > 0.
\end{equation}
Let $\sigma_g^2$ denote the pointwise variance of the underlying Gaussian field $G$. By prescribing the target mean $\mathbb{E}[\nu]=\nu_1$ and the target coefficient of variation $\mathrm{CoV}=\sigma_\nu/\nu_1$ for the physical viscosity, the corresponding Gaussian parameters are determined via standard lognormal relations \cite{ghanem1999nonlinear}:
\begin{equation}\label{eq:mu-sig-rel}
    \sigma_g^2=\ln\big(1+\mathrm{CoV}^2\big),\qquad\mu_\nu=\ln(\nu_1)-\tfrac{1}{2}\sigma_g^2.
\end{equation}

The spatial correlation of the underlying Gaussian field $G$ is governed by a separable exponential covariance kernel:
\begin{equation}\label{eq:cov-G-ours}
    \mathrm{Cov}\left(G(X_1),G(X_2)\right)=\sigma_g^2\exp\Big(-\tfrac{|x_2-x_1|}{L_x}-\tfrac{|y_2-y_1|}{L_y}\Big),\quad \text{for } X_1, X_2 \in D.
\end{equation}
The correlation lengths are chosen to be $25\%$ of the domain's global dimensions, yielding $L_x =  2$ and $L_y =  0.5$. Physically, these correlation lengths produce viscosity "patches" whose characteristic sizes are roughly one-quarter of the channel's length and height. Through equations \eqref{eq:nu-lognormal-ours} and \eqref{eq:cov-G-ours}, the instantaneous element Reynolds number, defined as $\mathrm{Re}(x;\bm\xi)=UL/\nu(x;\bm\xi)$, fluctuates around a nominal mean value of $\mathrm{Re}_1=100$, with the amplitude of these fluctuations strictly controlled by the specified $\mathrm{CoV}$. The underlying Gaussian field approximated using a truncated Karhunen-Loève (KL) expansion with $N=2$ terms. Subsequently, the nonlinear lognormal map \eqref{eq:nu-lognormal-ours} is projected onto a generalized polynomial chaos (gPC) basis $\{\psi_\ell(\bm\xi)\}$ composed of multivariate Hermite polynomials of total degree $p =3$. This yields the expansion:
\begin{equation}\label{eq:nu-gpc-ours}
    \nu(x;\bm\xi)\approx \sum_{\ell=1}^{n_\nu} \nu_\ell(x)\,\psi_\ell(\bm\xi), \qquad \text{where} \quad \nu_\ell(x)=\frac{\mathbb{E}\big[\nu(x;\bm\xi)\,\psi_\ell(\bm\xi)\big]}{\mathbb{E}\big[\psi_\ell^2(\bm\xi)\big]}.
\end{equation}
\noindent The subspace spanned by multivariate polynomials in $\{\xi_j\}_{j=1}^N$ of total degree $p=3$ yields dimension of $n_\xi = \binom{N+p}{p} = 10$~\cite{ghanem2003stochastic, xiu2002wiener}. The polynomial degree used for the expansion of the lognormal process is set to $2p=6$ to ensure a complete representation of the process in the discrete problem~\cite{ghanem1999nonlinear, sousedik2016stochastic}. Hence the stochastic Galerkin matrices $\bm{H}_{\ell}$, defined by $\bm{H}_{\ell}(i,j) = \mathbb{E}[\psi_i \psi_j \psi_\ell]$, are of size $10 \times 10$ with $\ell = 1, \ldots, 28$. A comprehensive summary of the physical, temporal, spatial, and stochastic parameters used in the Narrow Channel benchmark is provided in \cref{tab:setup-params}.
\begin{table}[ht]
\centering
\begin{tabular}{llc}
\toprule
\textbf{Category} & \textbf{Parameter Description} & \textbf{Value} \\
\midrule
\textit{Physical} 
& Final integration time ($t_{\mathrm{final}}$)          & $\SI{1.0}{\second}$ \\
\textit{~~\& Temporal}& Mean kinematic viscosity ($\nu_{1}$)                   & $0.01$ \\
& Viscosity coefficient of variation ($\mathrm{CoV}$)    & $10\%$ \\
& Number of temporal steps ($n_{t}$)                     & $40$ \\
\midrule
\textit{Spatial} 
& \# Velocity DOFs ($n_{u}$)                  & $1744$ \\
\textit{Discretization}& \# Pressure DOFs ($n_{p}$)                  & $281$ \\
& \# Total spatial DOFs ($n_x = n_u + n_p$) & $2025$ \\
\midrule
\textit{Stochastic} 
& Stochastic dimension ($N$)                             & $2$ \\
\textit{Discretization}& Solution gPC polynomial degree ($p$)                   & $3$ \\
& Viscosity expansion degree ($2p$)                      & $6$ \\
& Number of solution gPC modes ($n_{\xi}$)               & $10$ \\
& Number of viscosity expansion terms ($n_{\nu}$)        & $28$ \\
\midrule
\textit{Global System}
& \textbf{\# Total DOFs ($N_{\mathrm{total}}=n_t n_x n_\xi$)}  & \textbf{810,000} \\
\bottomrule
\end{tabular}
\caption{Summary of the physical parameters and the number of degrees of freedom (DOFs) used in the spatial, temporal, and stochastic discretizations for the Narrow Channel experiments.}
\label{tab:setup-params}
\end{table}

\subsection{Temporal Discretization and Sampling Methods}
\label{subsec:temporal-and-sampling}
\begin{figure}[ht]
   \centering \includegraphics[trim={.2cm 0 .2cm .1cm}, clip, width=0.9\textwidth]{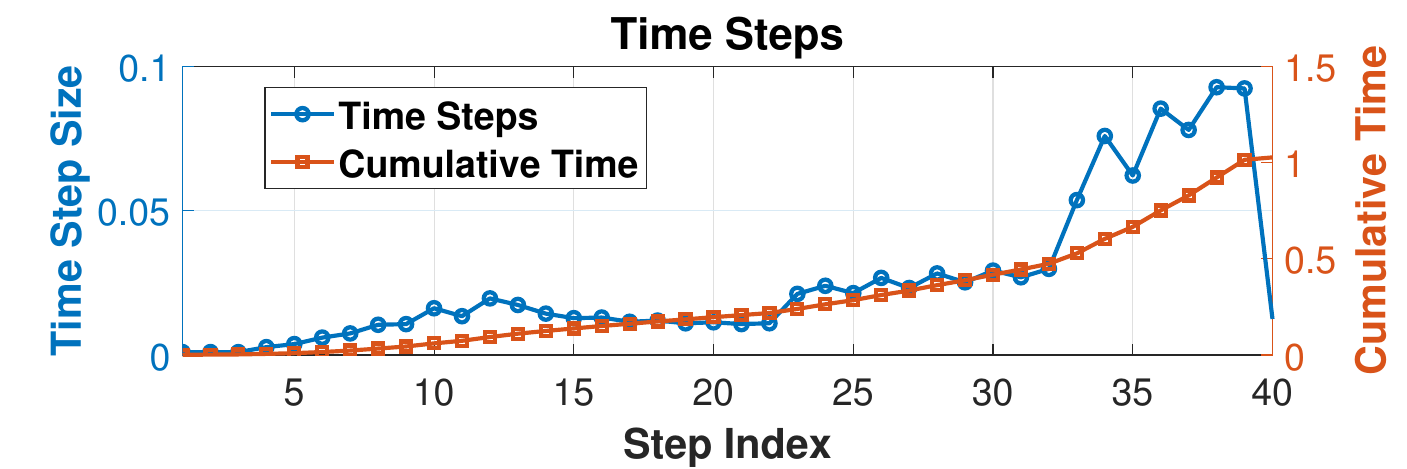}
  \caption{Temporal grids utilized for the integration of the Navier--Stokes equations. The secondary vertical axis denotes the cumulative simulation time.}
\label{fig:time-steps}
\end{figure}

 To efficiently capture the flow dynamics, particularly during the initial start-up phase, we first simulate the deterministic problem using the mean viscosity $\nu_1$. This deterministic run generates an adaptive sequence of time steps \cite{kay2010adaptive, sousedik2022stochastic}. \Cref{fig:time-steps} illustrates the sequence of time-step sizes determined by the deterministic solver, alongside the cumulative simulation time.

To validate the accuracy of the proposed stochastic Galerkin solver, we establish high-fidelity baselines using two widely adopted sampling techniques: the Monte Carlo (MC) method and the stochastic collocation (SC) method. Fundamentally, both sampling strategies rely on evaluating a series of completely decoupled, deterministic Navier--Stokes problems at specific sample points $\{\bm{\xi}^{(q)}\}$ in the random parameter space, yielding individual realizations of the flow field.

The \textbf{Monte Carlo} method generates a set of $N_{MC}$ sample points randomly in accordance with the underlying probability distribution of the random variables. The statistical moments of the flow field, such as the mean and variance, are then approximated using standard ensemble averaging over these deterministic realizations. In the experiments reported in this work, $N_{MC}=10^3$ samples are used to ensure a sufficiently accurate reference solution for comparison.

The \textbf{stochastic collocation} method constructs its samples using a structured, deterministic set of collocation points. These points are typically derived from sparse grid techniques to mitigate the curse of dimensionality \cite{gerstner1998numerical, novak1996high}. To enable a direct and consistent comparison between the SC baseline and our SGFEM formulation, we adopt a pseudospectral projection approach \cite{sousedik2016stochastic, xiu2010numerical}. In this framework, the corresponding gPC expansion coefficients for the velocity and pressure fields are explicitly computed via numerical quadrature. For example, the velocity coefficients are evaluated as
\begin{equation*}
  \textstyle u_{ik} = \sum_{q=1}^{N_q} \vec{u}^{(q)}(x_i)\psi_k(\bm{\xi}^{(q)})w^{(q)},
\end{equation*}
where $\vec{u}^{(q)}(x_i)$ is the deterministic solution at the $q$-th collocation point, and $w^{(q)}$ represents the corresponding quadrature weight. An analogous quadrature rule is applied to compute the pressure coefficients $p_{ik}$.

\subsection{Stochastic Flow Dynamics and Solver Convergence}
\label{subsec:validation-dynamics}
\begin{figure}[ht]
    \centering
    \includegraphics[trim={1.5cm .2cm 0.8cm .1cm}, clip, width=0.49\linewidth]{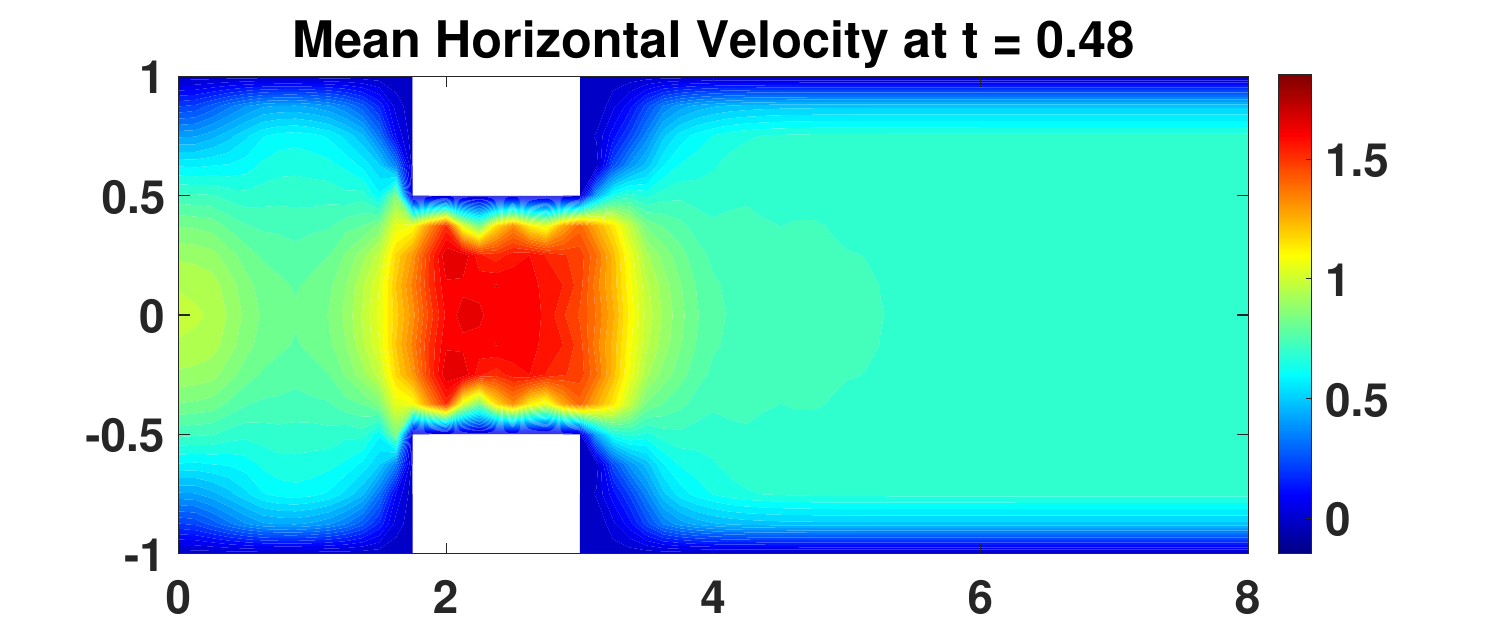}
    \includegraphics[trim={1.5cm .2cm 0.8cm .1cm}, clip, width=0.49\linewidth]{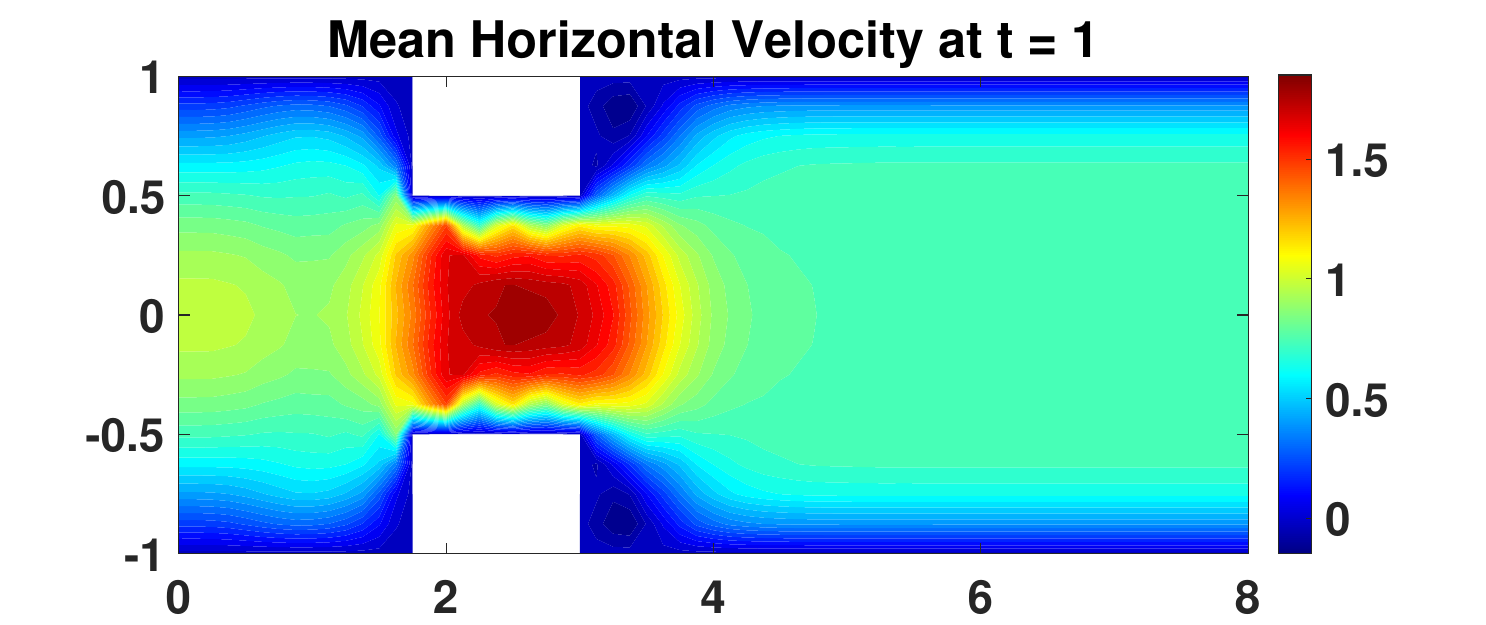} \\
    \includegraphics[trim={1.5cm .2cm 0.8cm .1cm}, clip, width=0.49\linewidth]{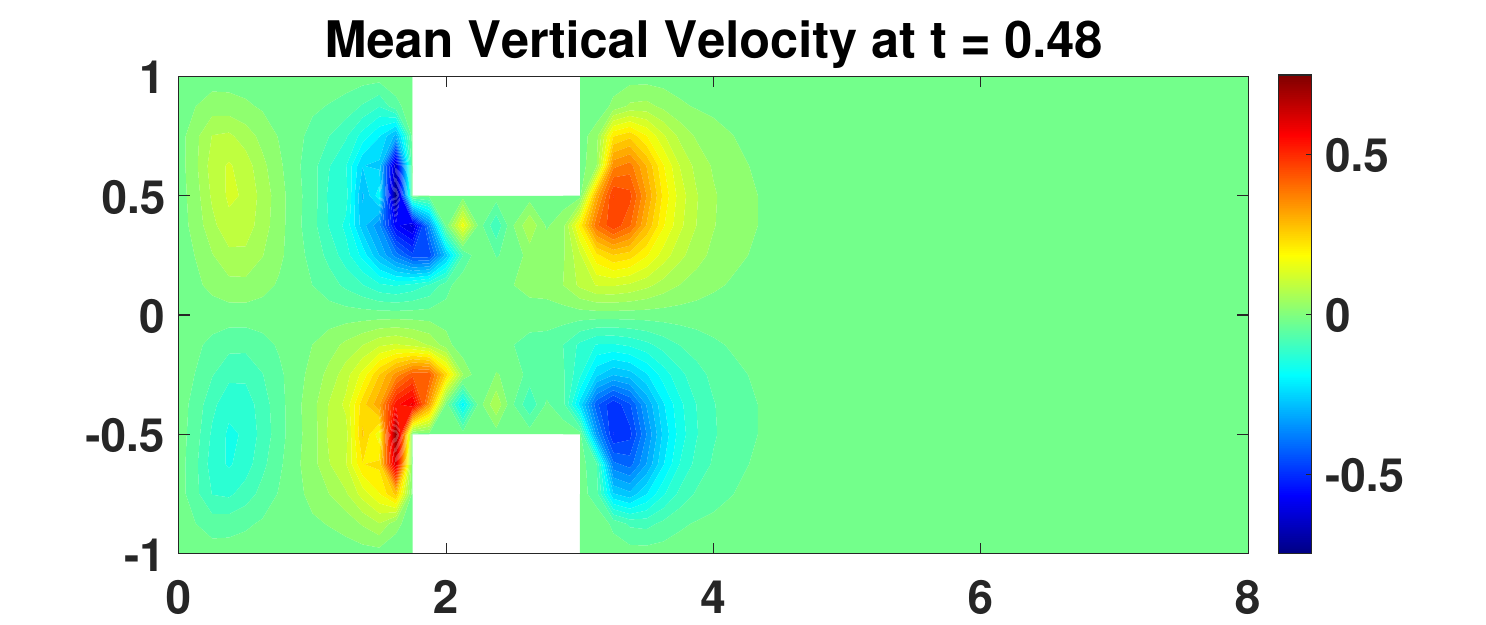}
    \includegraphics[trim={1.5cm .2cm 0.8cm .1cm}, clip, width=0.49\linewidth]{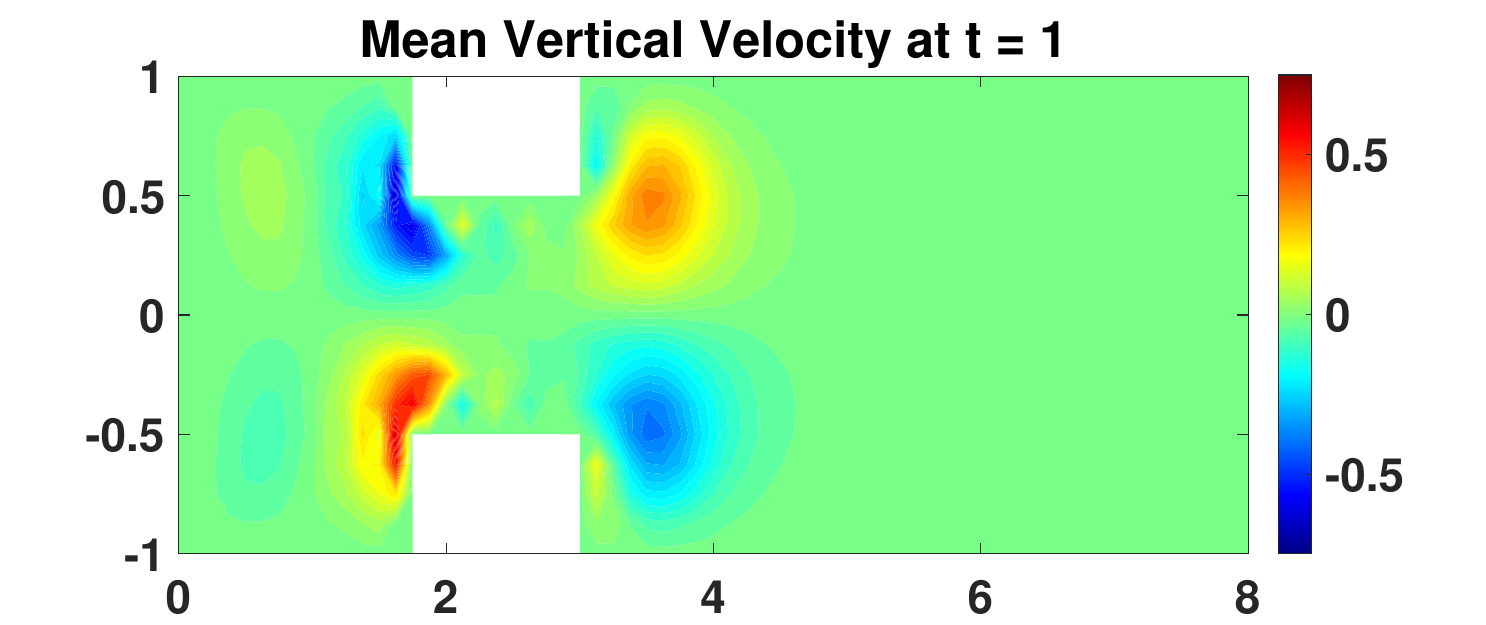} \\
    \includegraphics[trim={1.5cm .2cm 0.8cm .1cm}, clip, width=0.49\linewidth]{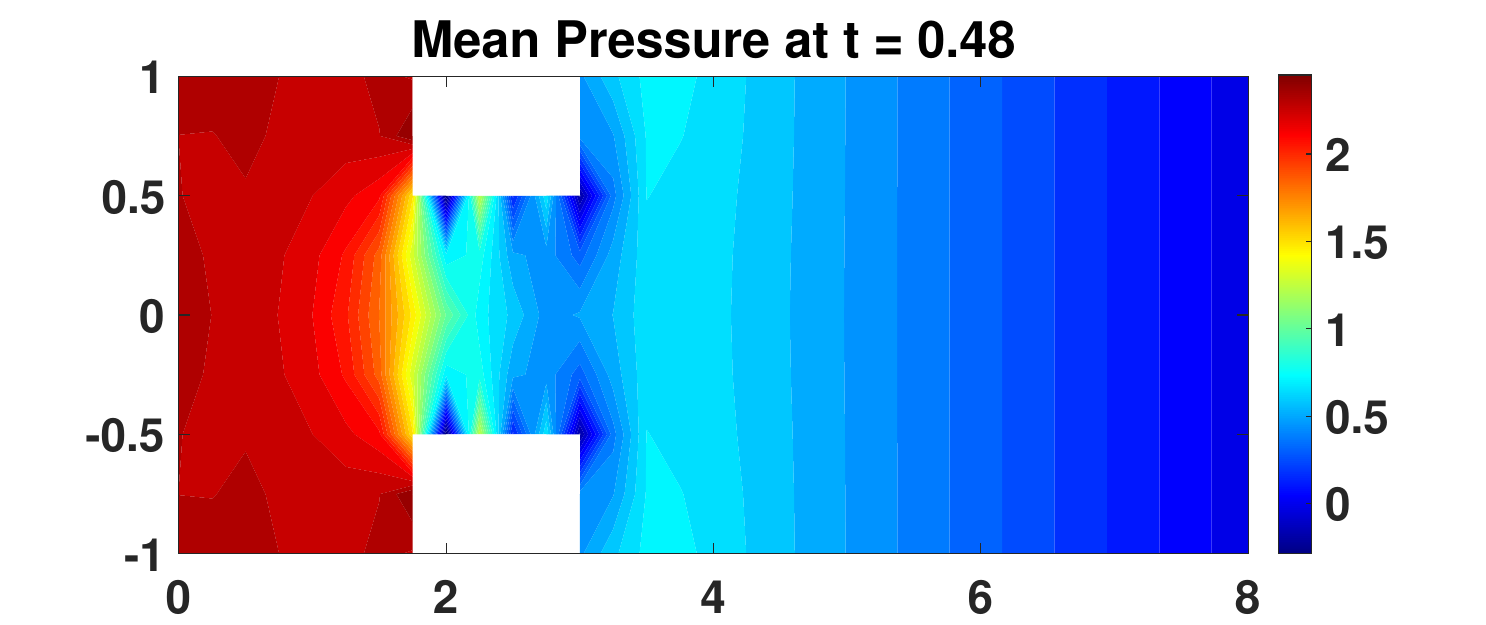}
    \includegraphics[trim={1.5cm .2cm 0.8cm .1cm}, clip, width=0.49\linewidth]{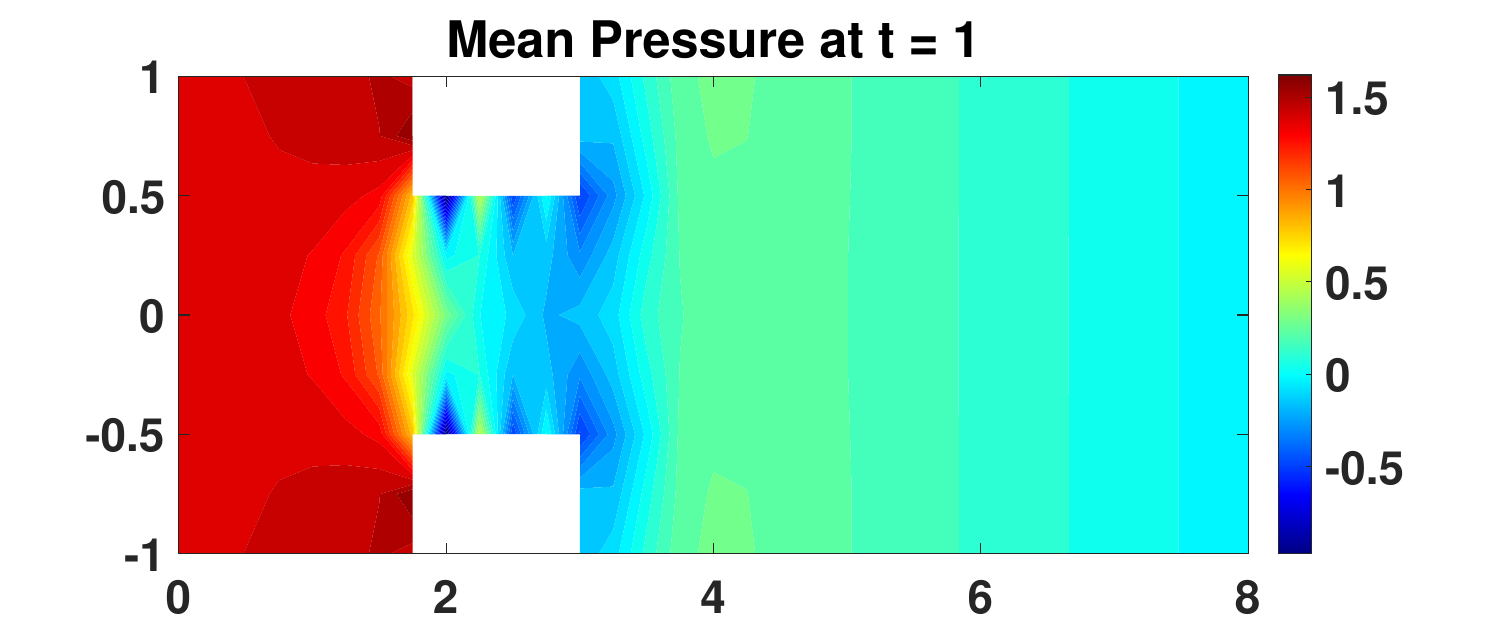}
    \caption{Spatial contour plots of the mean velocity and pressure fields. The left column displays the intermediate time $t=\SI{0.48}{\second}$, and the right column displays the final time $t=\SI{1.0}{\second}$. From top to bottom:  horizontal velocity, vertical velocity, and pressure.}
    \label{fig:mean-fields}
\end{figure}

\begin{figure}[H]
    \centering
    \includegraphics[trim={1.5cm .2cm 0.8cm .1cm}, clip, width=0.49\linewidth]{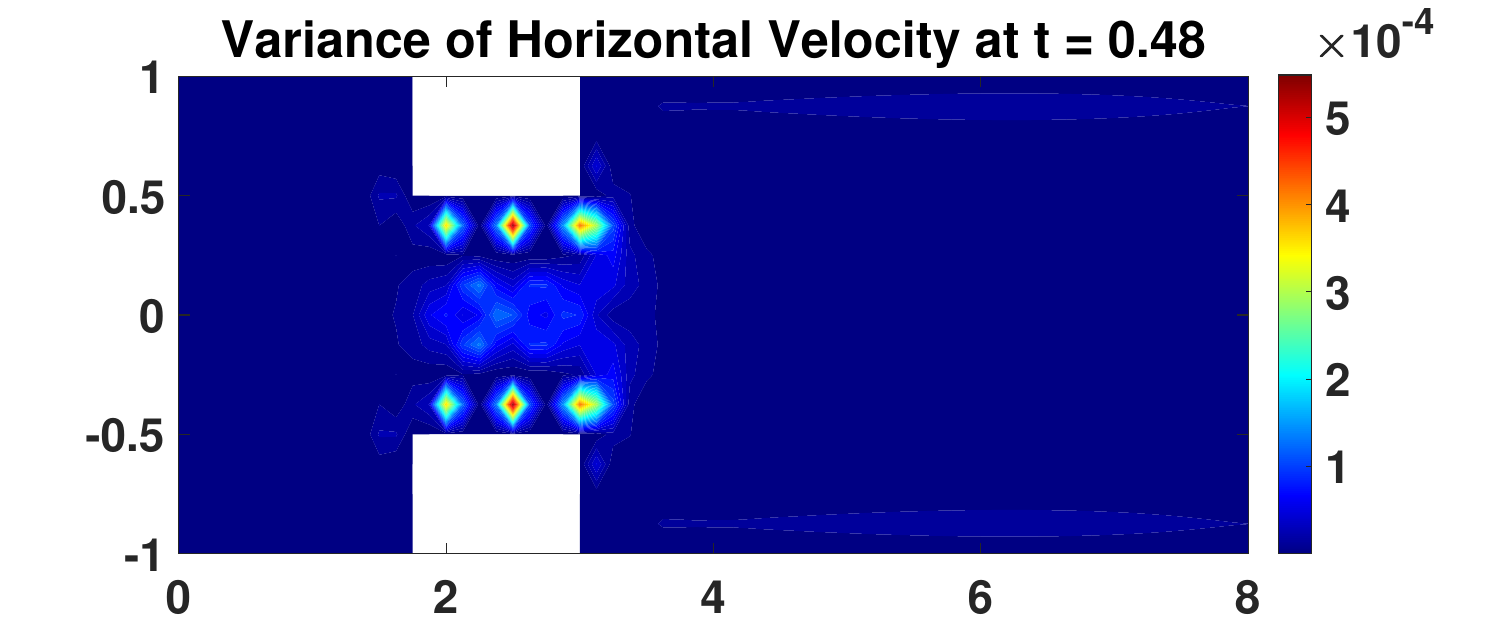}
    \includegraphics[trim={1.5cm .2cm 0.8cm .1cm}, clip, width=0.49\linewidth]{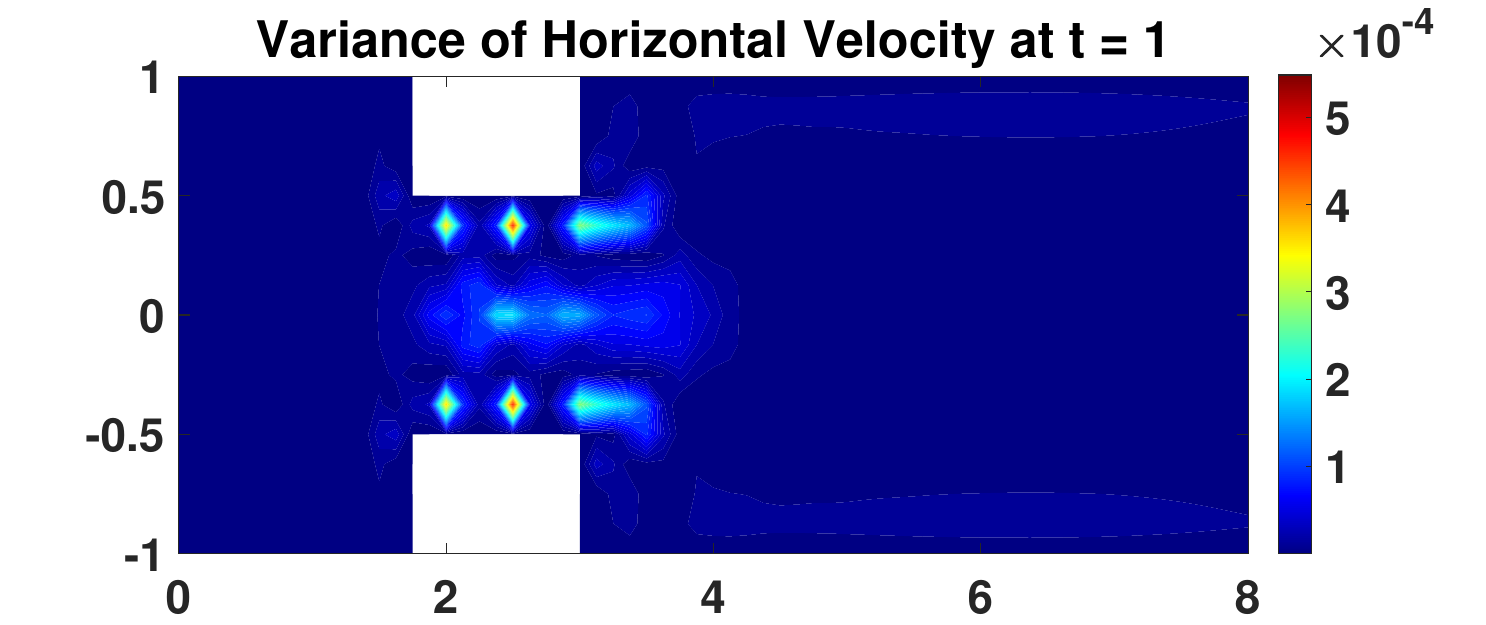} \\
    \includegraphics[trim={1.5cm .2cm 0.8cm .1cm}, clip, width=0.49\linewidth]{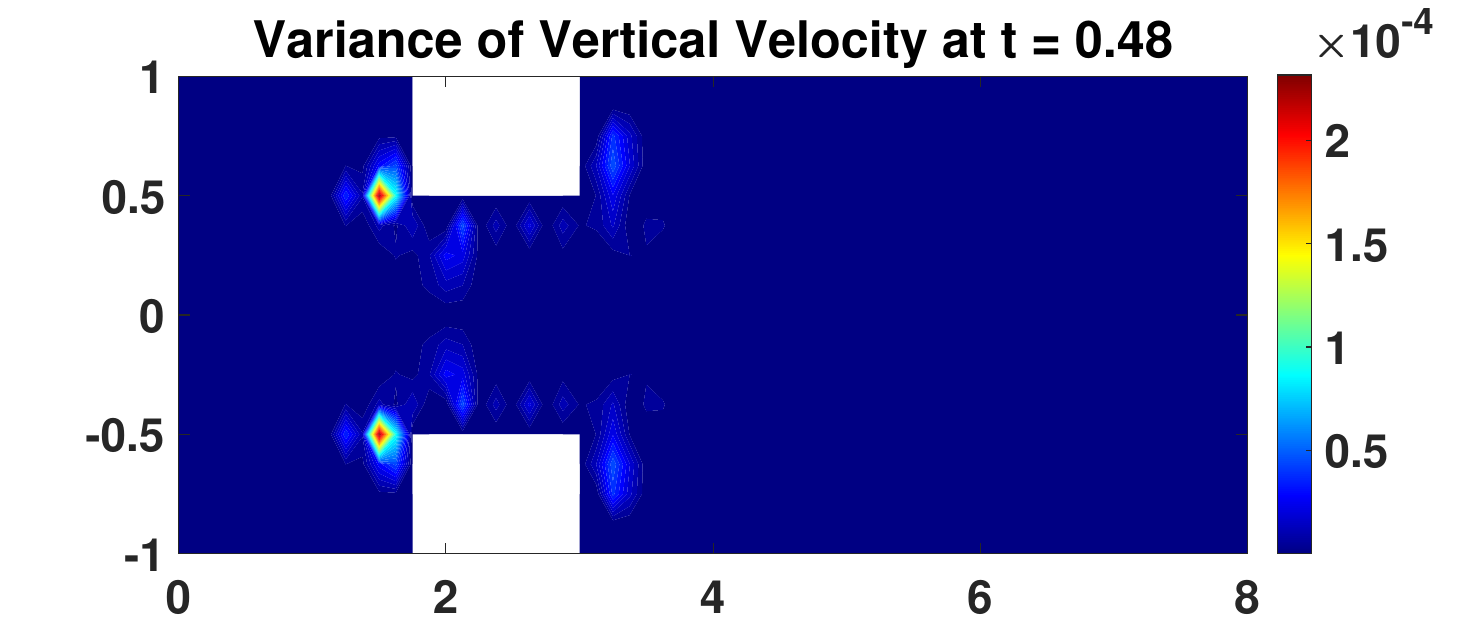}
    \includegraphics[trim={1.5cm .2cm 0.8cm .1cm}, clip, width=0.49\linewidth]{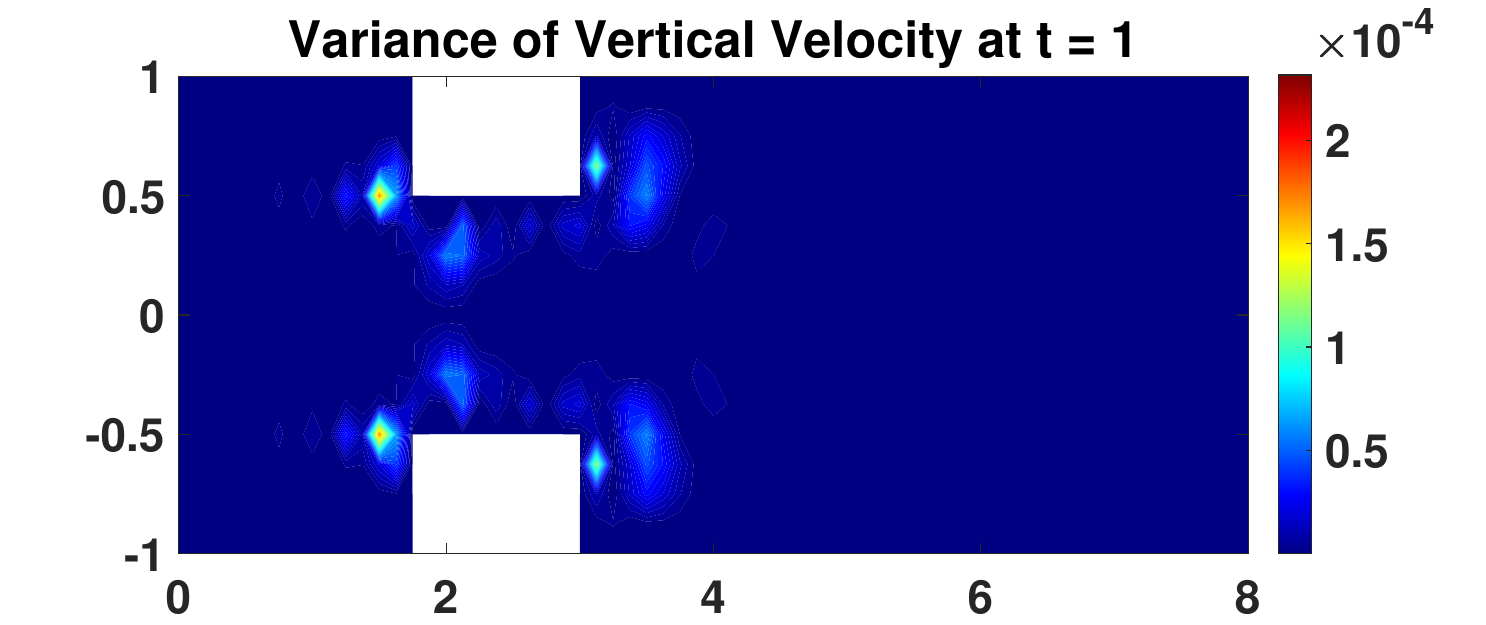} \\
    \includegraphics[trim={1.5cm .2cm 0.8cm .1cm}, clip, width=0.49\linewidth]{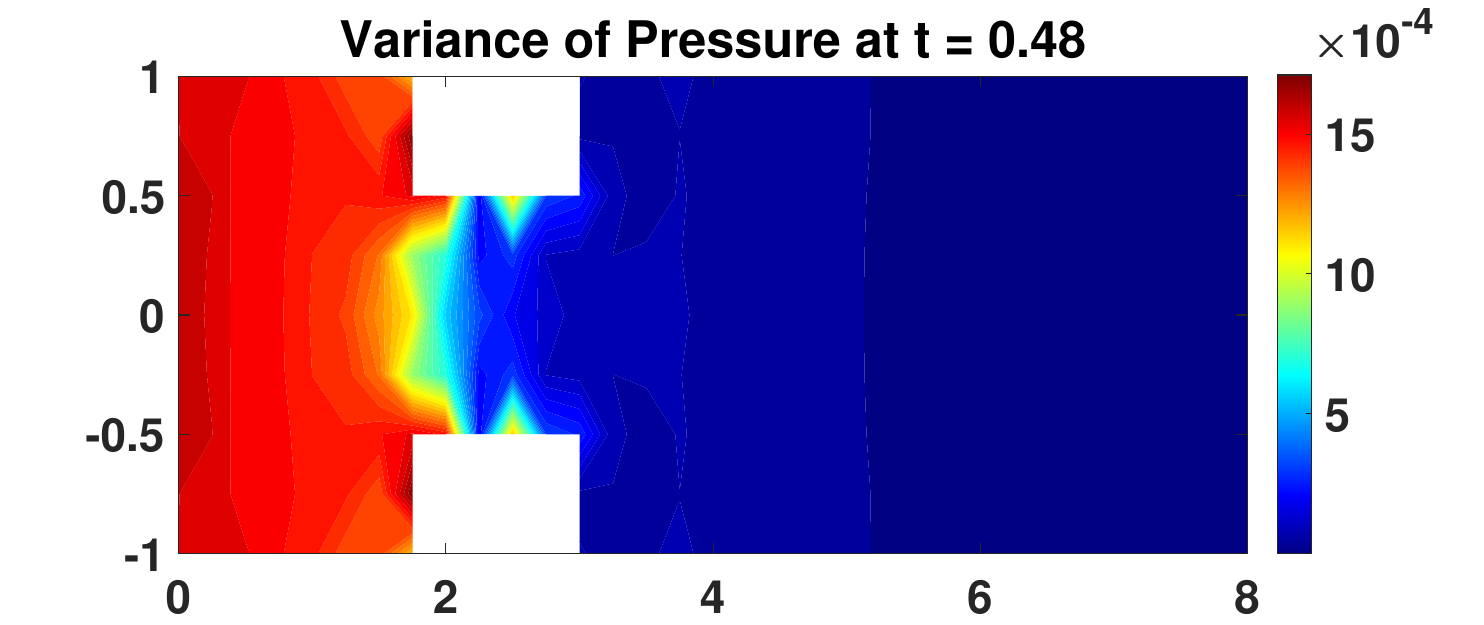}
    \includegraphics[trim={1.5cm .2cm 0.8cm .1cm}, clip, width=0.49\linewidth]{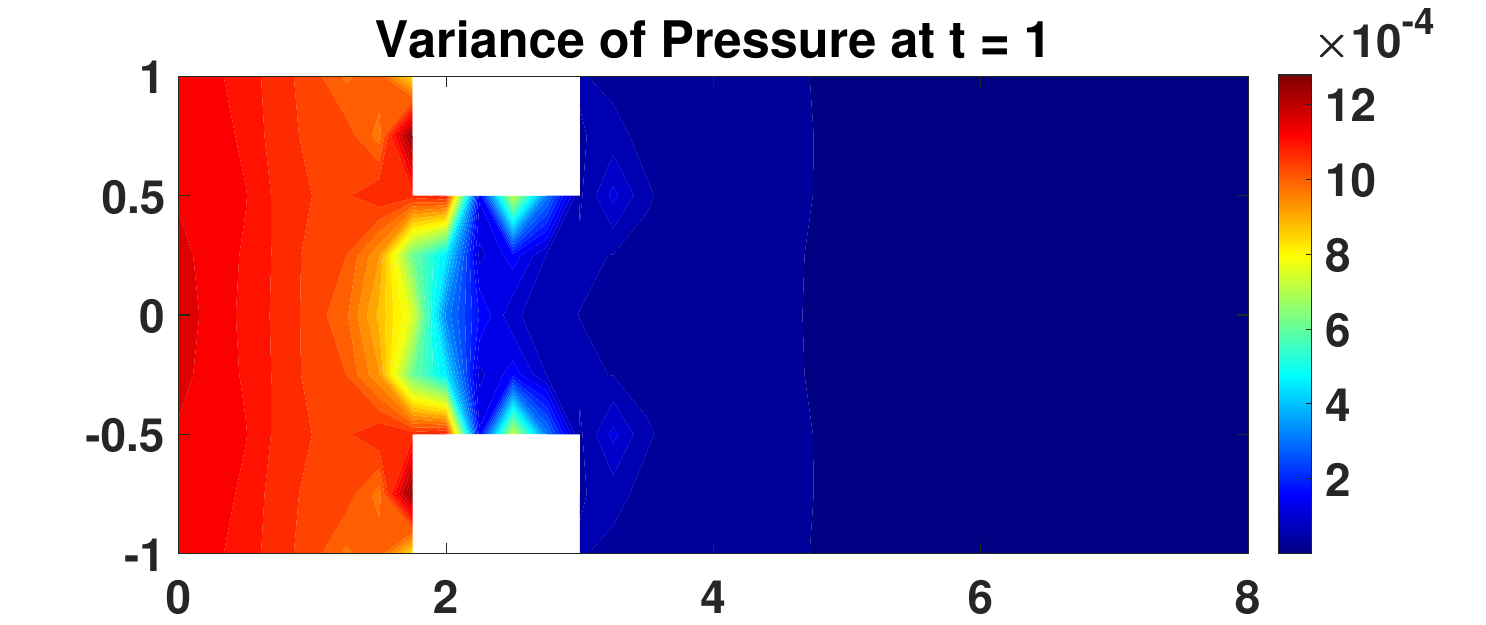}
    \caption{Spatial contour plots illustrating the variance fields. The left column displays the intermediate time $t=\SI{0.48}{\second}$, and the right column displays the final time $t=\SI{1.0}{\second}$. From top to bottom: horizontal velocity, vertical velocity, and pressure.}
    \label{fig:variance-fields}
\end{figure}

\begin{figure}[ht]
    \centering
    \includegraphics[trim={1cm .2cm 2cm .2cm}, clip, width=0.495\linewidth]{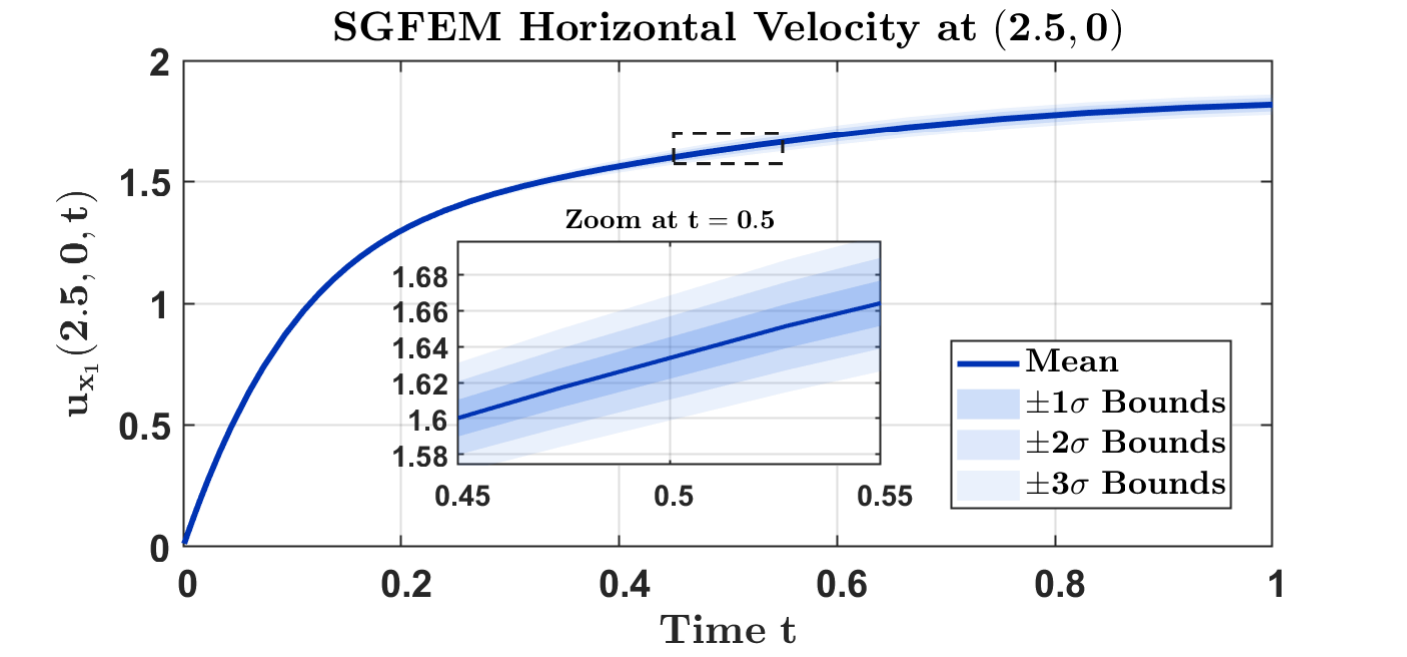}
    \hfill
    \includegraphics[trim={1cm .2cm 2cm .2cm}, clip, width=0.495\linewidth]{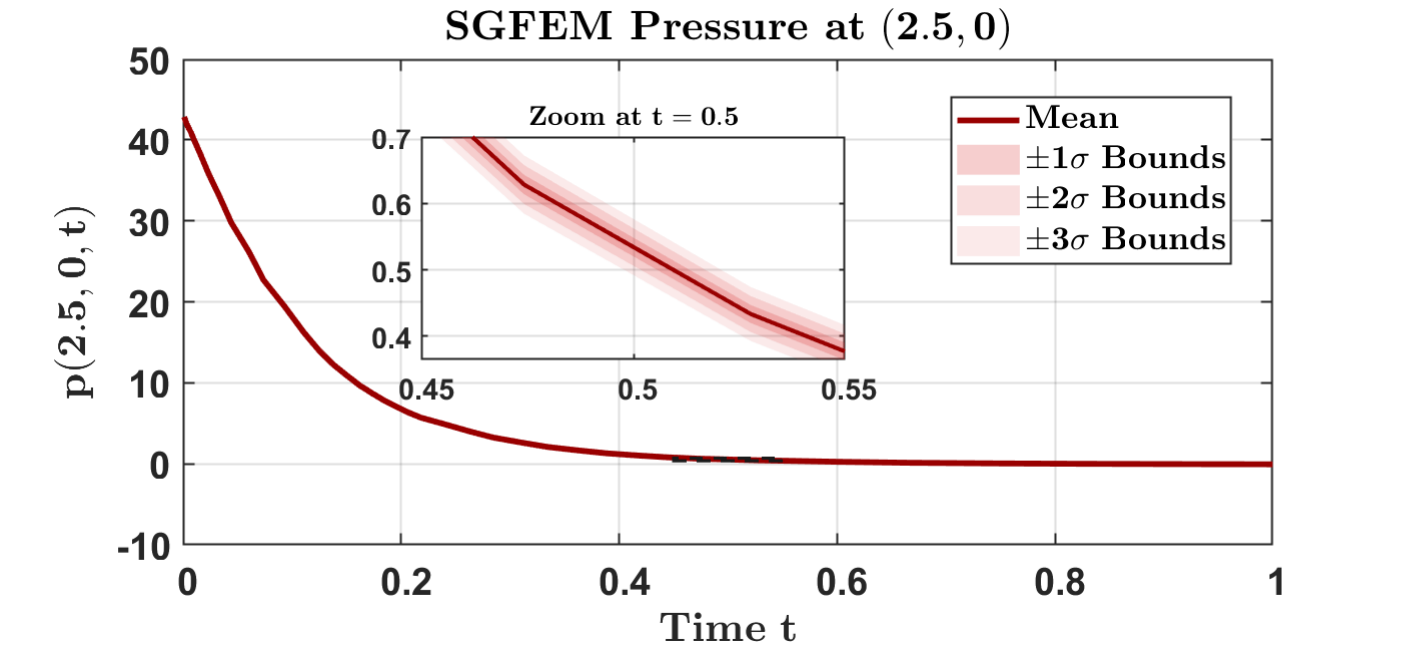}
    \caption{Evolution of the horizontal velocity (left) and pressure (right) evaluated continuously from $t=\SI{0}{\second}$ to $t=\SI{1.0}{\second}$ at the spatial coordinate $(2.5, 0)$. The shaded regions represent the standard deviation bounds ($\pm \sigma$) around the mean trajectory.}
    \label{fig:overtime-point}
\end{figure}

The mean and variance of the velocity and pressure fields are shown in \cref{fig:mean-fields} and \cref{fig:variance-fields}, respectively, at two time snapshots: an intermediate time $t=\SI{0.48}{\second}$, corresponding to the time step closest to $t=\SI{0.5}{\second}$, and the final time $t=\SI{1.0}{\second}$. As the parabolic inflow drives the flow, the mean horizontal velocity reaches its largest values in the narrow constriction. To illustrate the temporal evolution of this uncertainty, \cref{fig:overtime-point} tracks the horizontal velocity and pressure continuously from $t=\SI{0}{\second}$ to $t=\SI{1.0}{\second}$ at the $(2.5, 0)$, located near the middle of the constriction. The shaded regions represent the standard deviation bounds ($\pm \sigma$). \Cref{tab:error-matrix} reports the relative $L_2$ errors of the mean fields obtained by the three methods considered.

\begin{table}[H]
    \centering
    \caption{Relative $L_2$ errors ($\|\mathbb{E}_i - \mathbb{E}_j\| / \|\mathbb{E}_j\|$) for the mean velocity vector field $\bm{u}$ and pressure field $p$, comparing the SGFEM, SC, and MC methodologies.}
    \label{tab:error-matrix}
    \begin{tabular}{l@{\hskip 0.2in}cc@{\hskip 0.2in}cc}
    \toprule
    & \multicolumn{2}{c@{\hskip 0.2in}}{\textbf{Velocity Field ($\bm{u}$)}} & \multicolumn{2}{c}{\textbf{Pressure Field ($p$)}} \\
    \cmidrule(r{0.5in}){2-3} \cmidrule{4-5}
    \textbf{Method}  & \textbf{SC} & \textbf{MC}  & \textbf{SC} & \textbf{MC} \\
\midrule
      \textbf{SGFEM}   & $1.12 \times 10^{-4}$ & $1.44 \times 10^{-4}$  & $1.24 \times 10^{-05}$ & $3.49 \times 10^{-5}$ \\
      \textbf{SC}     & & $1.11 \times 10^{-4}$ & & $3.57 \times 10^{-5}$ \\
      \bottomrule
    \end{tabular}
\end{table}

\begin{figure}[ht]
    \centering
    \includegraphics[trim={.6cm .2cm .6cm .2cm}, clip,width=0.243\textwidth]{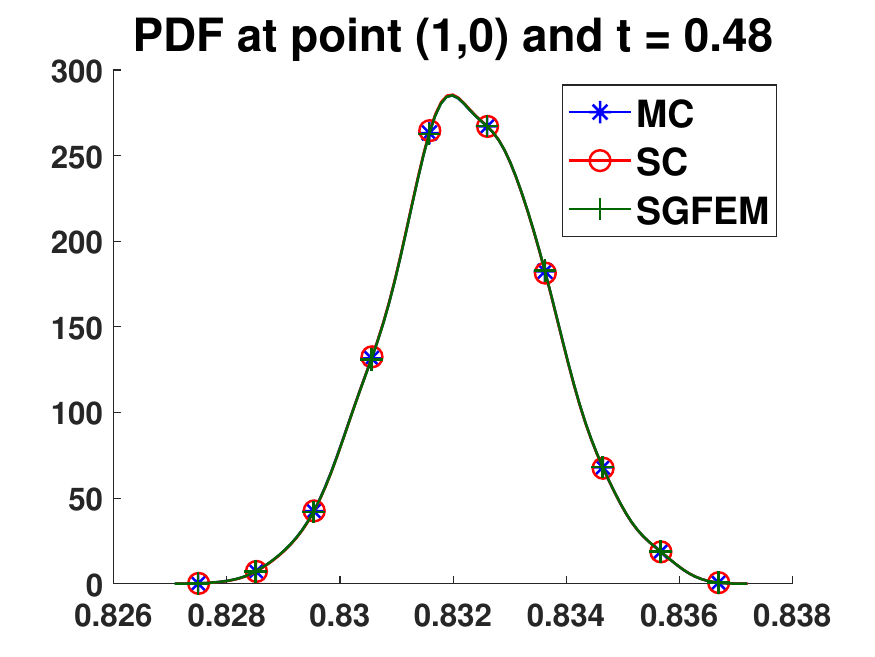}
    \includegraphics[trim={.6cm .2cm .6cm .2cm}, clip,width=0.243\textwidth]{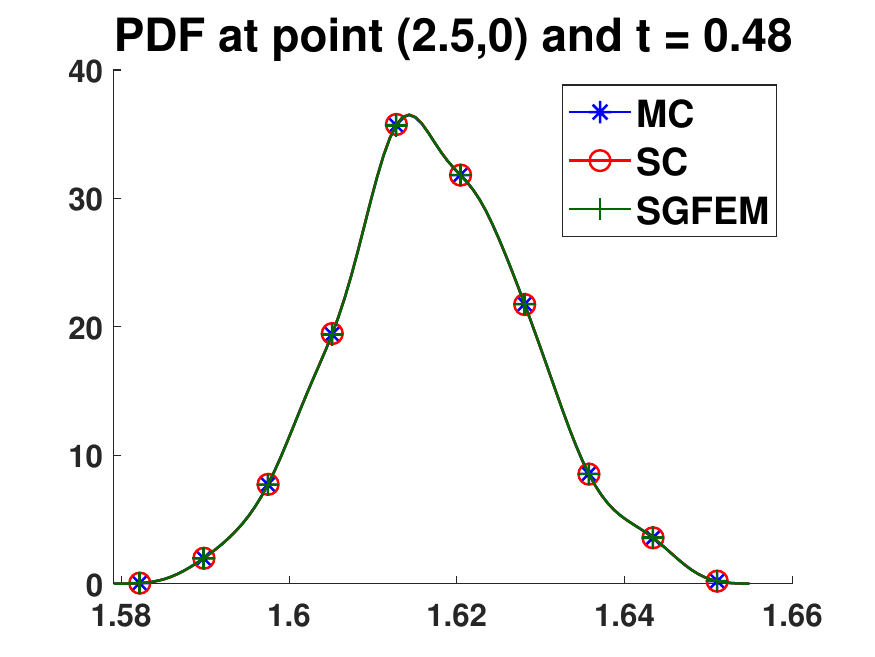}
    \includegraphics[trim={.6cm .2cm .6cm .2cm}, clip,width=0.243\textwidth]{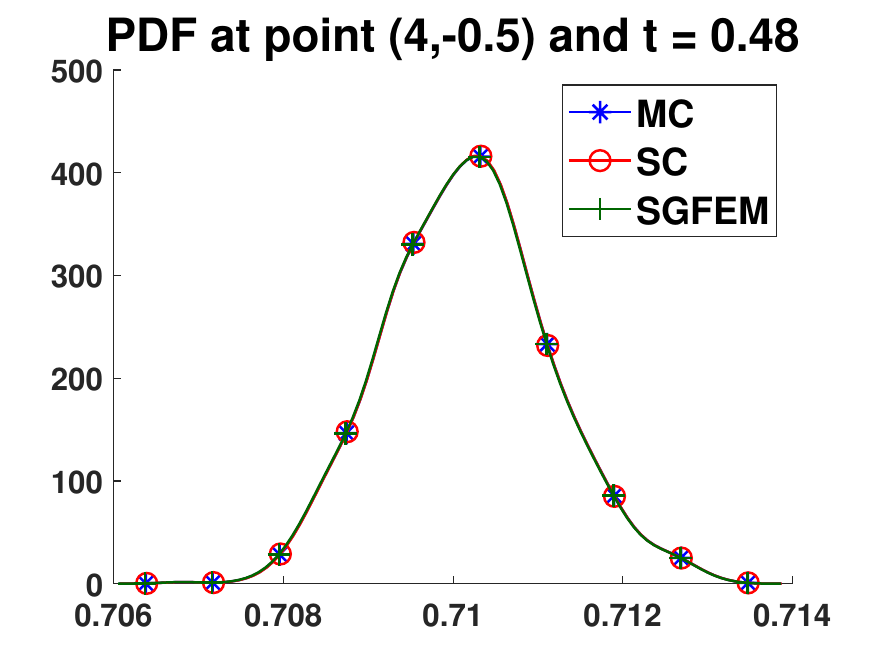}
    \includegraphics[trim={.6cm .2cm .6cm .2cm}, clip,width=0.243\textwidth]{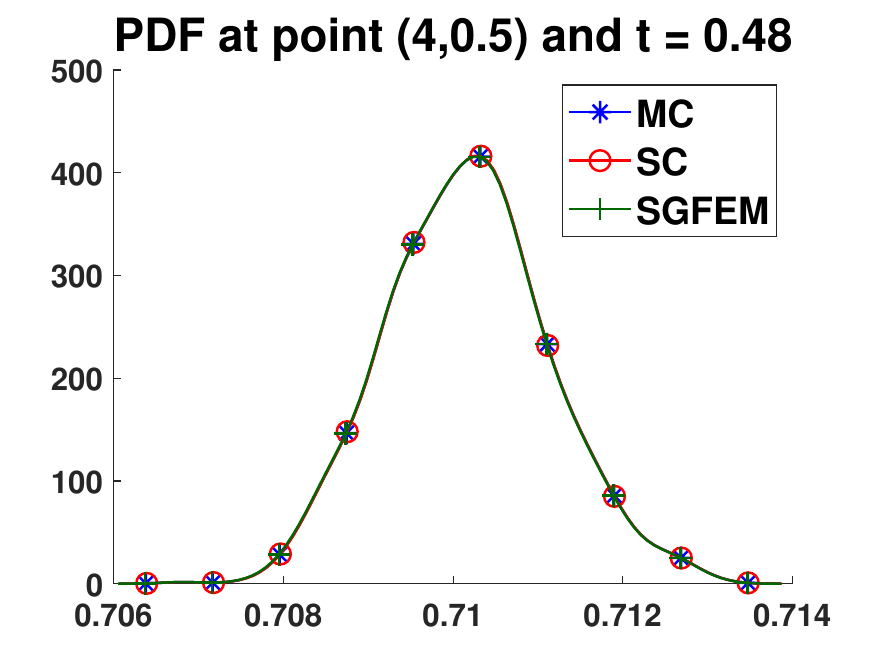}\\
    \vspace{0.4cm}
    \includegraphics[trim={.6cm .2cm .6cm .2cm}, clip,width=0.243\textwidth]{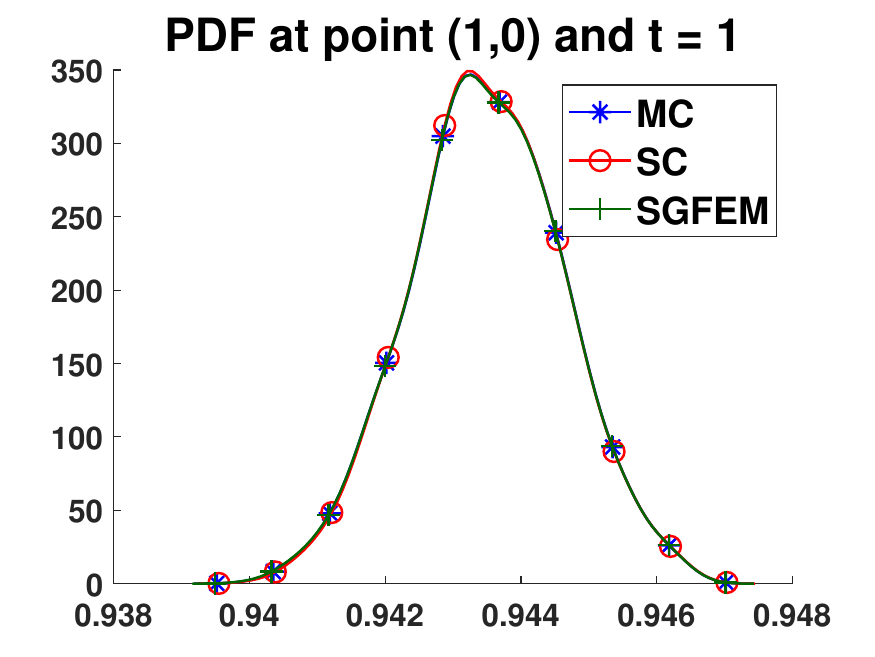}
    \includegraphics[trim={.6cm .2cm .6cm .2cm}, clip,width=0.243\textwidth]{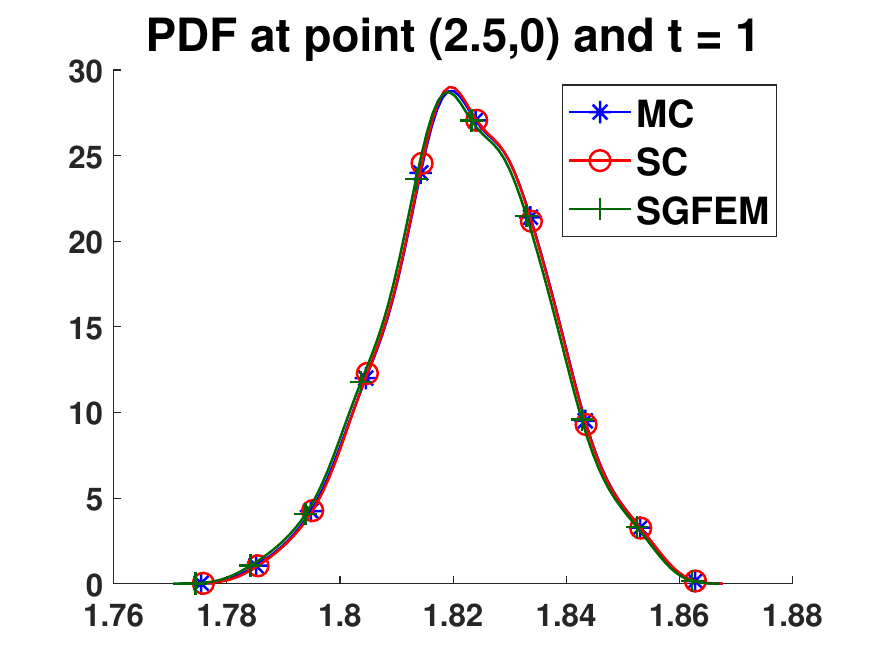}
    \includegraphics[trim={.6cm .2cm .6cm .2cm}, clip,width=0.243\textwidth]{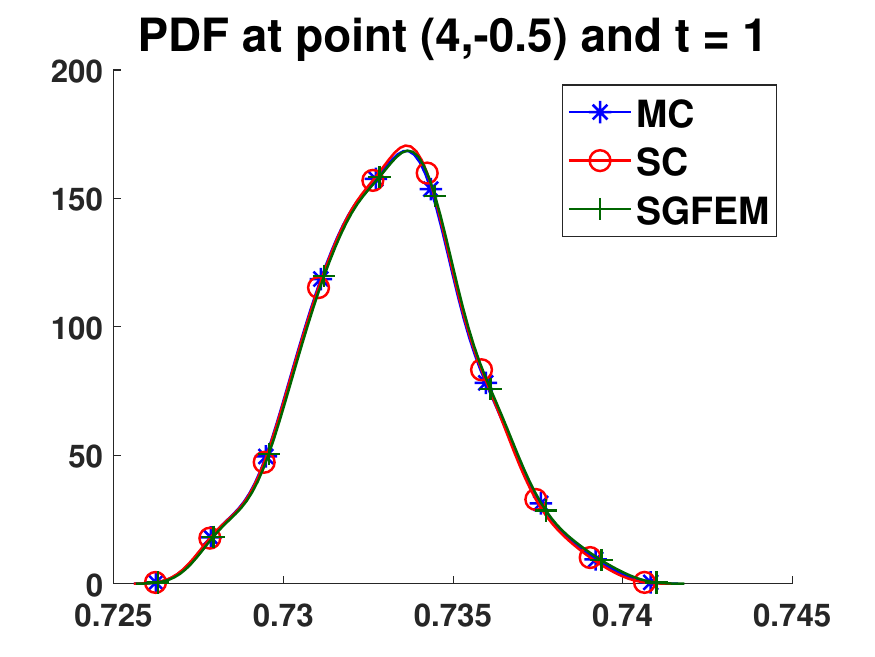}
    \includegraphics[trim={.6cm .2cm .6cm .2cm}, clip,width=0.243\textwidth]{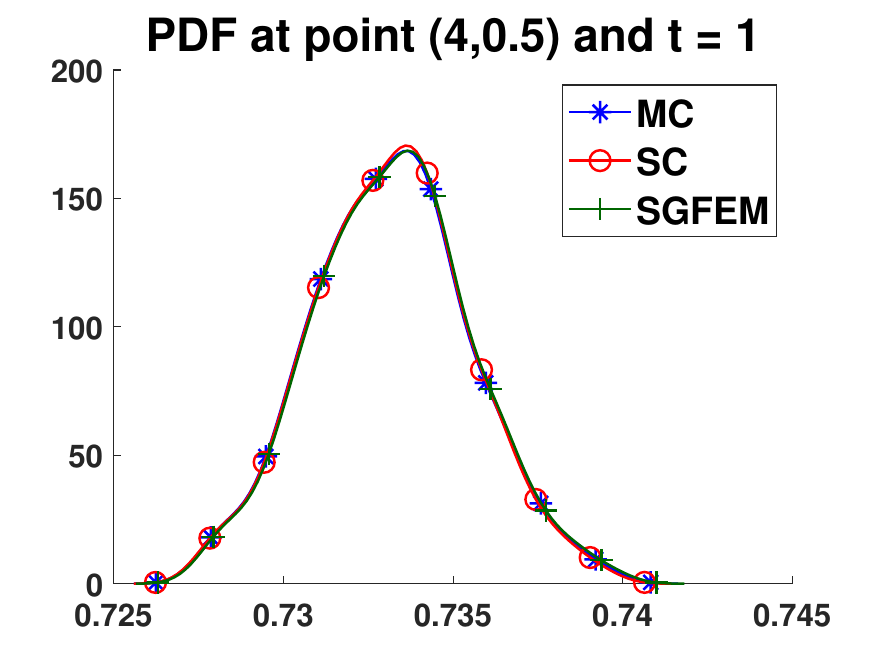}
    \caption{Probability Density Functions (PDFs) evaluated at four specific spatial points: $(1,0)$, $(2.5,0)$, $(4, -0.5)$, and $(4,0.5)$. The top and bottom rows display distributions at $t=\SI{0.48}{\second}$ and $t=\SI{1.0}{\second}$, respectively.}
    \label{fig:pdf-comparisons}
\end{figure}

\begin{figure}[ht]
    \centering
    \includegraphics[trim={.4cm 0 1.7cm .2cm}, clip, width=0.51\linewidth]{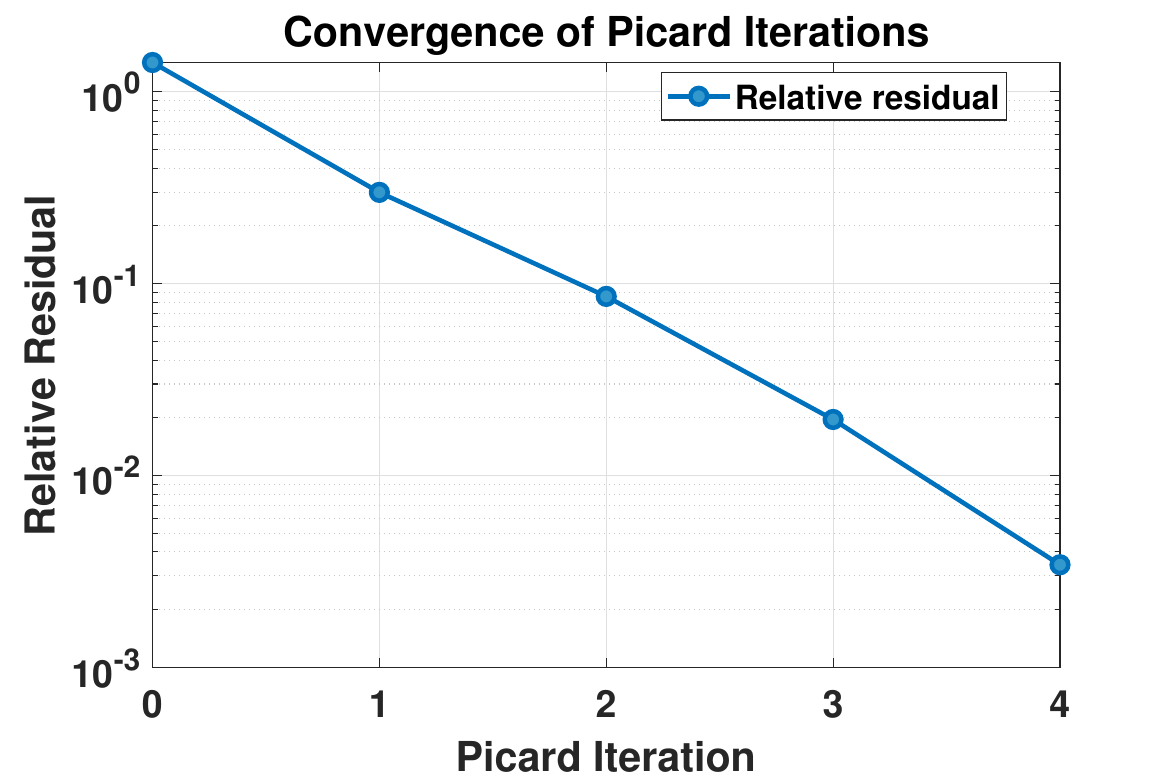}
    \includegraphics[trim={1.2cm 0 1.55cm .2cm}, clip, width=0.48\linewidth]{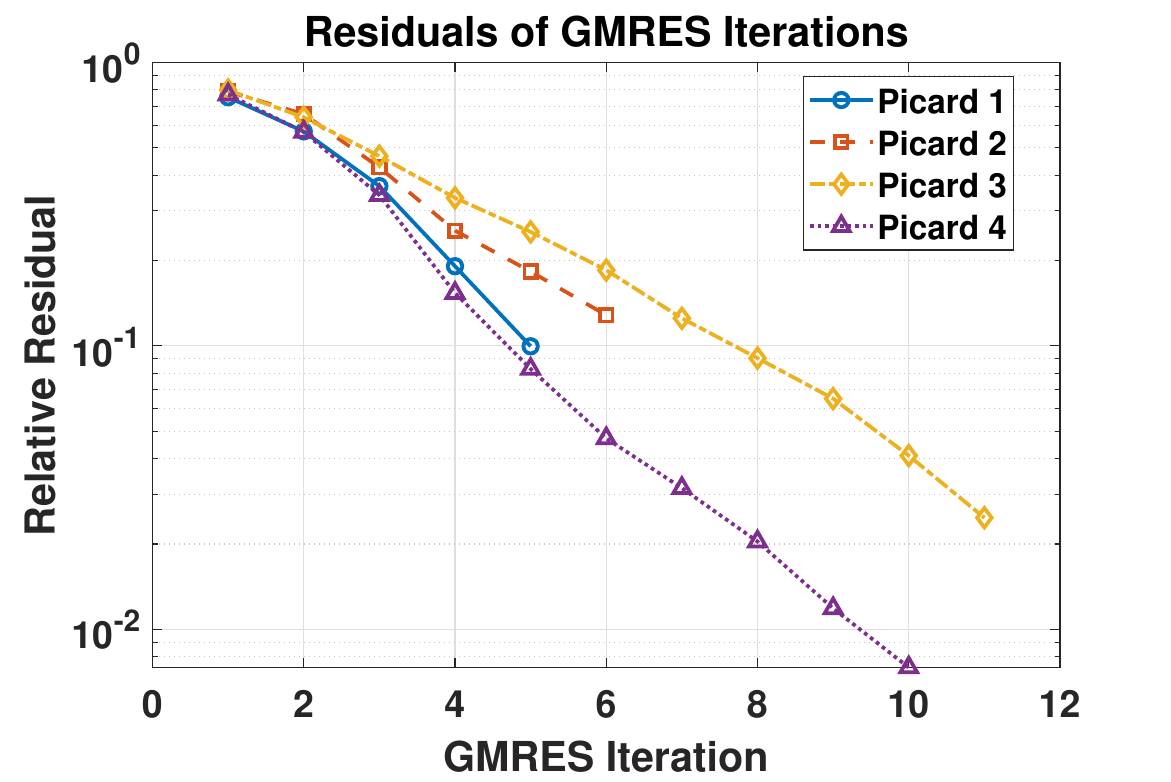}
    \caption{Left: Convergence history of the nonlinear Picard iteration, displaying the relative residual on a semi-logarithmic scale. Right: Convergence history of the outer GMRES iterations corresponding to each distinct Picard step. The outer GMRES successfully drives the linear residual down without stagnation.}
    \label{fig:convergence-behavior}
\end{figure}

Probability density functions (PDFs) show the method's ability to capture local stochastic behaviors. The PDFs of the horizontal velocity at four specific spatial coordinates of interest shown in  \Cref{fig:pdf-comparisons}. In \cref{fig:convergence-behavior} the left panel shows the relative nonlinear residual, illustrating the linear convergence of the Picard iteration to the prescribed tolerance. The right panel shows the outer GMRES residuals at each Picard step..

\subsection{Low-Rank Approximations and Tolerances}
\label{subsec:algorithmic-performance}
\begin{figure}[ht]
    \centering
        \includegraphics[trim={0.7cm 0cm 2.5cm 0}, clip, width=0.47\textwidth]{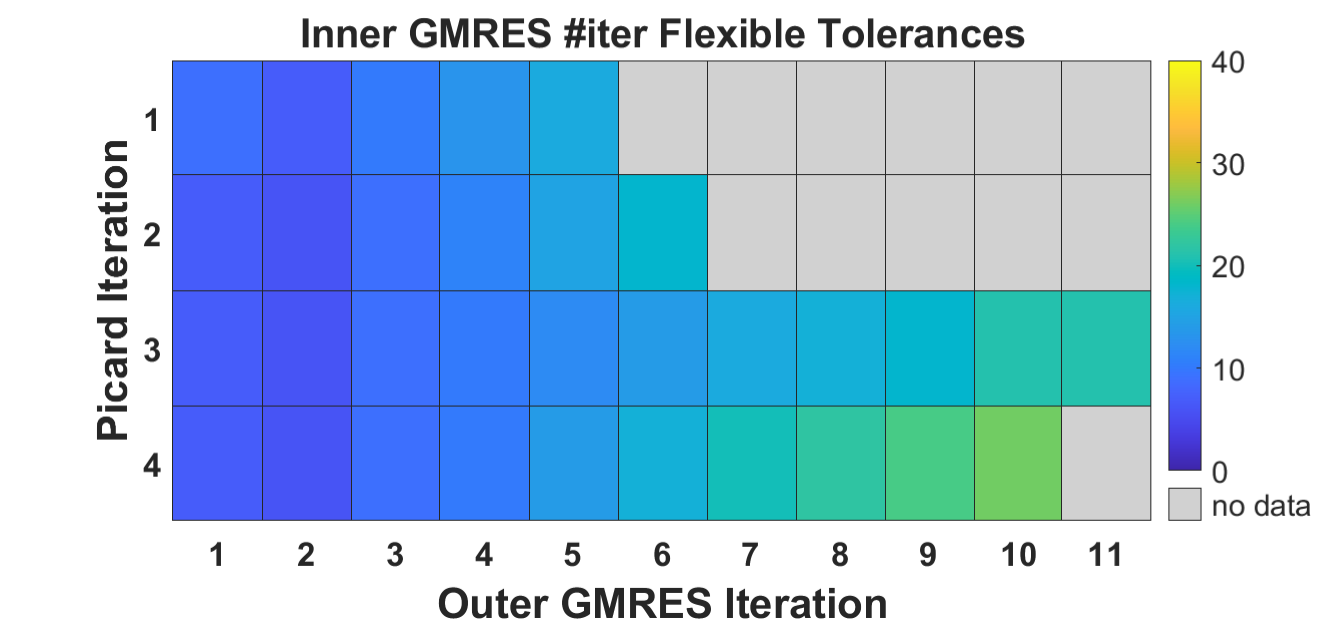}
        \hfill
        \includegraphics[trim={1.1cm 0 0 0}, clip, width=0.52\textwidth]{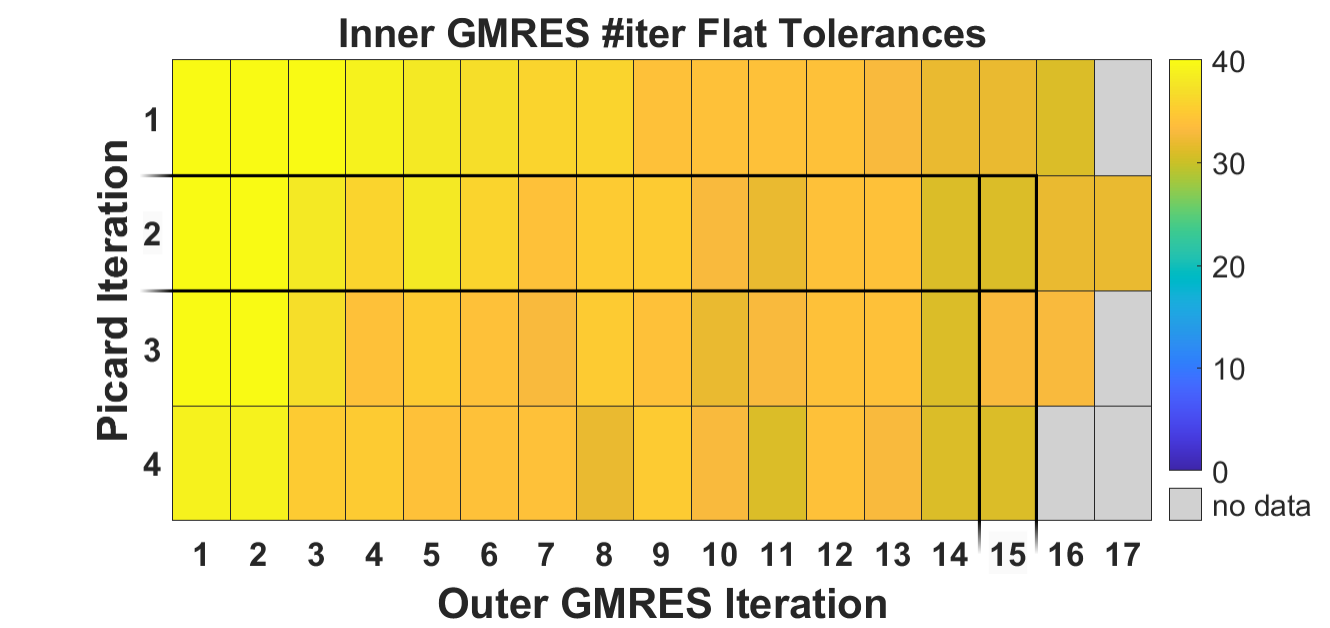}
    \caption{Heatmaps detailing the number of inner GMRES iterations nested within each outer GMRES step across the sequential Picard iterations. Darker cells represent fewer inner iterations. The flexible (Left) and flat (right) tolerance strategies are shown. The flexible strategy dramatically reduces the computational burden during early nonlinear updates.}
    \label{fig:iteration-counts}
\end{figure}

Using fixed tolerances for all inner solves and tensor truncations can lead to unnecessary work during the early Picard steps. For this reason, we consider a \textit{flexible} tolerance strategy, described in~\Cref{alg:picard_krylov}, in which the tolerances for the inner GMRES iterations and TT-rounding are scaled by the current residual. We compare this strategy with a \textit{flat} tolerance approach, where the same tolerances are imposed at every iteration. In particular, the flexible tolerances are chosen so that, once the outer residual reaches its target level, the corresponding inner tolerances match their flat counterparts. In \cref{tab:tt-params_tolerances}, the specific tolerance values used in both strategies are listed. \Cref{fig:iteration-counts} shows side-by-side heatmaps of the number of inner GMRES iterations required within each outer GMRES step, across the Picard iterations. The flat tolerance approach requires a consistently high number of iterations from the start, needlessly expending computational effort. The flexible strategy gives $20$-fold speed-up in computational time.

\begin{table}[ht]
    \centering
    \caption{Comparison of tolerances for the flexible (relative) and flat (fixed) strategies.}
    \label{tab:tt-params_tolerances}
    \begin{tabular}{lcccc}
    \toprule
    \textbf{Solver level} & \multicolumn{2}{c}{\textbf{Iteration Tolerance}} & \multicolumn{2}{c}{\textbf{TT-Rounding Tolerance}} \\
    \cmidrule(lr){2-3} \cmidrule(lr){4-5}
    & \textbf{Flex} & \textbf{Flat} & \textbf{Flex} & \textbf{Flat} \\
    \midrule
    Picard & $1\times10^{-2}$ & $1\times10^{-2}$ & $5\times10^{-5}$ & $5\times10^{-7}$ \\
    Outer GMRES & $1\times10^{-1}$ & $1\times10^{-3}$ & $5\times10^{-4}$ & $5\times10^{-7}$ \\
    Inner GMRES & $5\times10^{-1}$ & $5\times10^{-4}$ & $1\times10^{-3}$ & $5\times10^{-7}$ \\
    \bottomrule
    \end{tabular}
\end{table}

\begin{figure}[ht]
    \centering
        \includegraphics[trim={0 0 .7cm 0}, clip, width=0.51\textwidth]{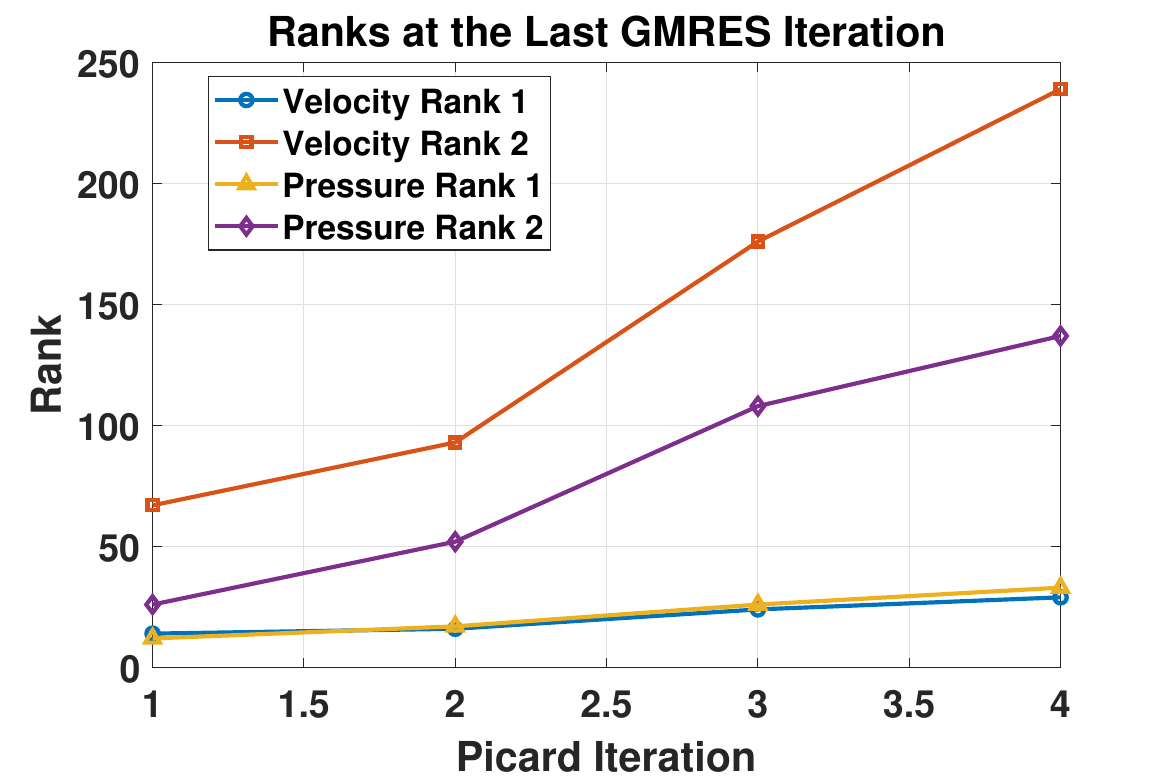}
        \includegraphics[trim={1.05cm 0 .7cm 0}, clip, width=0.48\textwidth]{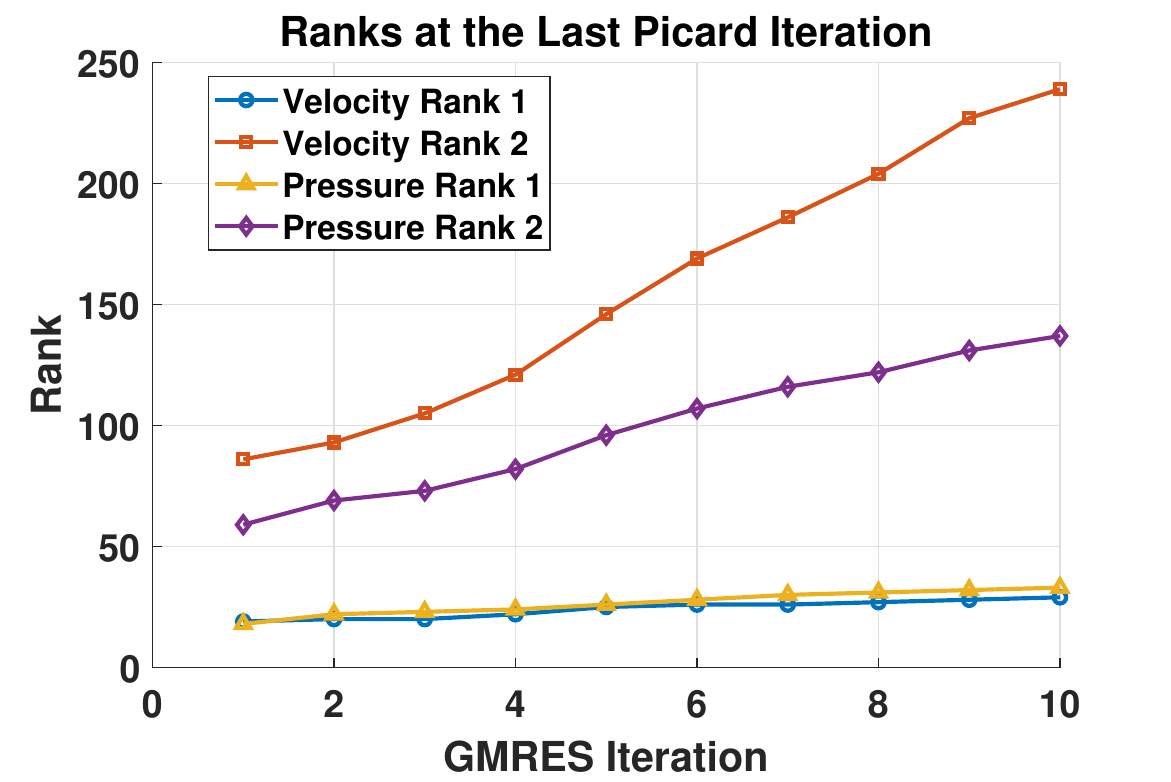}
    \caption{Evolution of the TT ranks ($\kappa_1, \kappa_2$). The top row illustrates the flexible tolerance strategy, where ranks begin at highly compressed levels and dynamically grow only as accuracy demands increase. The bottom row illustrates the flat tolerance strategy, where ranks remain uniformly bloated from the initial iterations.}
    \label{fig:ranks-evolution}
\end{figure}

A similar relative tolerance strategy is applied to the TT-SVD (\cref{alg:TT_SVD}) rounding operations to dynamically control the ranks of the solution tensor ($\kappa_1, \kappa_2$) and prevent memory bloat. \Cref{fig:ranks-evolution} illustrates the evolution of the TT ranks during the last outer GMRES (left) and last Picard iteration (right) for the solve strategy with flexible tolerances. The growth in the ranks reflects the tightening of the relative rounding tolerances; as the nonlinear residual becomes smaller, the TT-SVD truncation threshold decreases.

\begin{figure}[ht]
    \centering
    \includegraphics[trim={2.4cm 0 2.4cm .2cm}, clip, width=\linewidth]{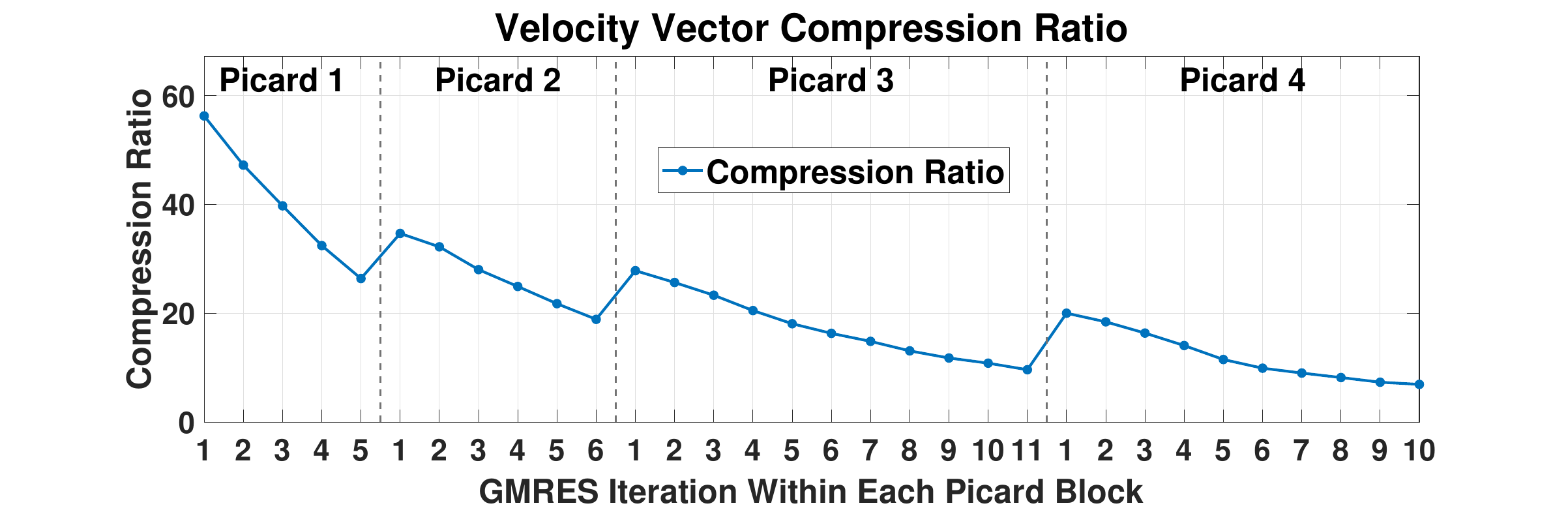}
    \caption{Evolution of the velocity vector compression ratio across the inner GMRES iterations, grouped by overarching Picard steps. The flexible tolerance strategy induces dynamic compression ratios due to the relative rounding thresholds. TT format offers at least compression ratio of $7:1$, reducing the memory usage by $86\%$.}
    \label{fig:compression-ratio}
\end{figure}

The memory savings achieved by the TT format is computed by the compression ratio of the velocity solution vector throughout the nested iterative solve. The compression ratio is defined as the ratio between the full size vector and the size of the TT representation. Let $\kappa_1$ and $\kappa_2$ denote the ranks of the velocity tensor. The compression ratio is calculated as:
\begin{equation}
    \text{Compression Ratio} = \frac{\text{size of full vector} }{
        \text{size of TT representation}
    } = \frac{n_t n_\xi n_u}{\kappa_1 n_t + \kappa_1 \kappa_2 n_\xi + \kappa_2 n_u}.
\end{equation}
\Cref{fig:compression-ratio} shows the compression ratio evaluated at each outer GMRES iteration, segmented by the corresponding Picard iterations. The nonlinearity and stochastic couplings inherently introduce finer structural details into the flow field, causing the TT ranks ($\kappa_1, \kappa_2$) to grow. The periodic increases at the beginning of each new Picard block is because the tolerance is scaled relative to the current residual. Most importantly, even at the very last GMRES iteration where the solution achieves its maximum required complexity, the TT representation maintains a compression ratio of $7:1$. This corresponds to an $86\%$ reduction in overall memory usage.

\subsection{Comparative Analysis of Preconditioners}
\label{subsec:precon_comparison}
The overall efficiency of the monolithic solver heavily relies on the convergence of the inner GMRES iterations that apply $\mathbb{F}_u^{-1}$. This inversion is computationally demanding due to the intricate coupling across temporal, spatial, and stochastic dimensions. We compare four distinct preconditioning strategies, ranging from simple mean-based approximations to our proposed structure-preserving low-rank formulations:

\begin{itemize}
    \item The \textit{Mass-Matrix Preconditioner} assumes that the mass matrix with time coupling dominates the system, while the nonlinear convection term and stochastic coupling are neglected.
    \begin{equation}
    \mathbb{P}_{\text{mass}}^{-1} =  \bm{T}^{-1} \otimes \bm{I}_{n_\xi} \otimes \bm{M}_*^{-1}.
    \label{eq:prec-simple}
    \end{equation}

    \item The \textit{Mean--Field Preconditioner}, \eqref{eq:innerGMRES_Elman}, introduced in \cite{elman2020low}.

    \item The \textit{Stochastic Neumann-Series Preconditioner}, defined in \eqref{eq:prec-neumann-expansion}, uses CP approximations of the operators together with a damped Neumann expansion. In the experiments, the damping parameter is chosen heuristically as $\delta=0.1$, and the CP rank is set to $3$.

    \item The \textit{CP1 Preconditioner}, defined in \eqref{eq:prec-cp1-inv}, approximates the global operator $\mathbb{F}_u$ using a rank-one CP decomposition.
\end{itemize}

\begin{table}[htbp]
\centering
\caption{Total computational time (CPU seconds) across various inner preconditioners and time-stepping regimes. CP1 preconditioner yields the most efficient overall performance, while the Mean-Field approach fails under nonuniform stepping. The Mean-Field preconditioner algorithm was terminated prematurely due to non-convergence within $100$ iterations under the nonuniform time-stepping regime.}
\label{tab:comp_times}
\begin{tabular}{@{}lccccc|cc@{}}
\toprule
& \multicolumn{5}{c}{Coefficient of Variation (CoV) $\sigma_\nu = 0.1$} & \multicolumn{2}{c}{$\sigma_\nu = 0.01$} \\ \cmidrule(lr){2-8}
& \textbf{Adaptive} & \multicolumn{6}{c}{\textbf{Uniform Time Steps $n_t$}} \\ \cmidrule(lr){3-8} 
\textbf{Preconditioner} & \textbf{(Non-unif.)} & $2^6$ & $2^7$ & $2^8$ & $2^9$ &$2^7$ & $2^9$ \\ \midrule
Mass Matrix \eqref{eq:prec-simple}   & ~2440 & 2776 & 3573 & 4314 & 4683&1833&2348 \\
Mean Field \eqref{eq:innerGMRES_Elman}   & $> 25240$ & 1318 & 2284 & 4141 & 5620 &915&1901\\
Neumann Series \eqref{eq:prec-neumann-expansion} & ~1375 & 2453 & 2584 & 3013 & 4303 &1320&1471\\
CP1 \eqref{eq:prec-cp1-inv}        & ~1629 & 1234 & 1322 & 1868 & 2030 &641&971\\ \bottomrule
\end{tabular}
\end{table}

\begin{figure}[ht]
    \centering
    \includegraphics[trim={0 0.9cm 2.5cm 0}, clip, width=0.46\textwidth]{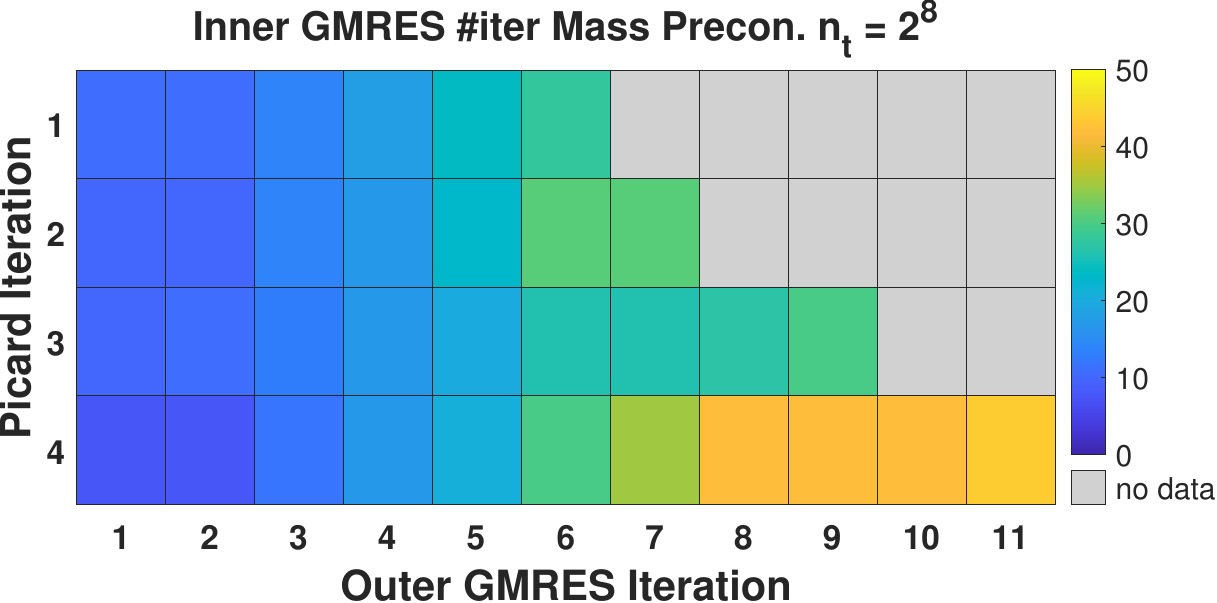}
    \hfill
    \includegraphics[trim={0.8cm 0.9cm 0 0}, clip, width=0.5\textwidth]{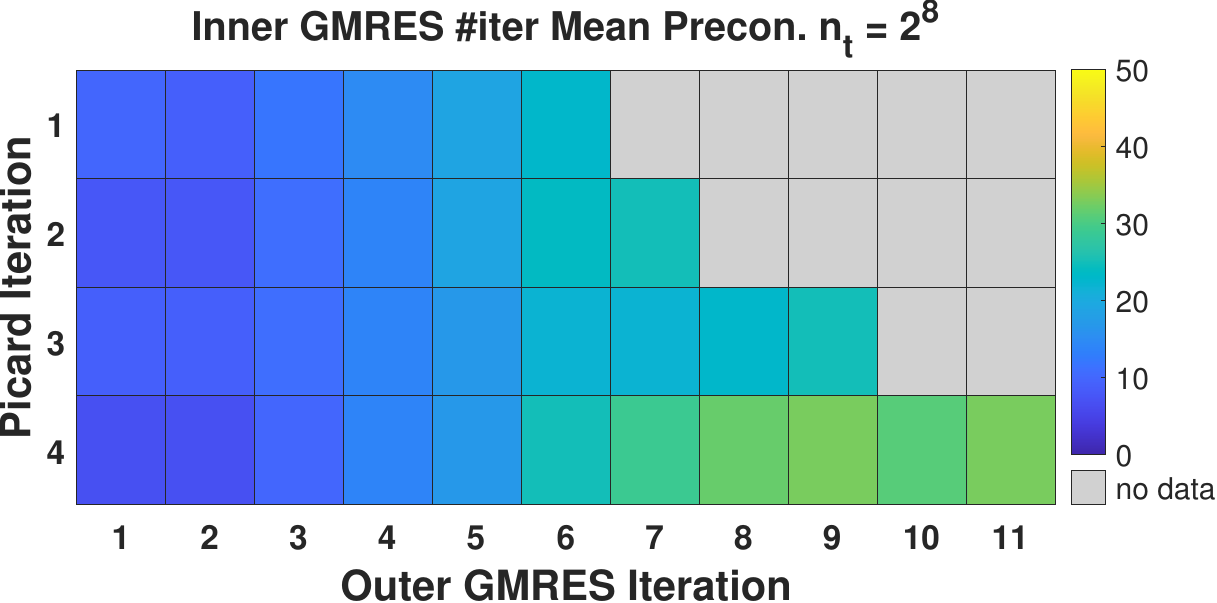}
   \\
   \vspace{0.2cm}
    \includegraphics[trim={0 0 2.5cm 0}, clip, width=0.46\textwidth]{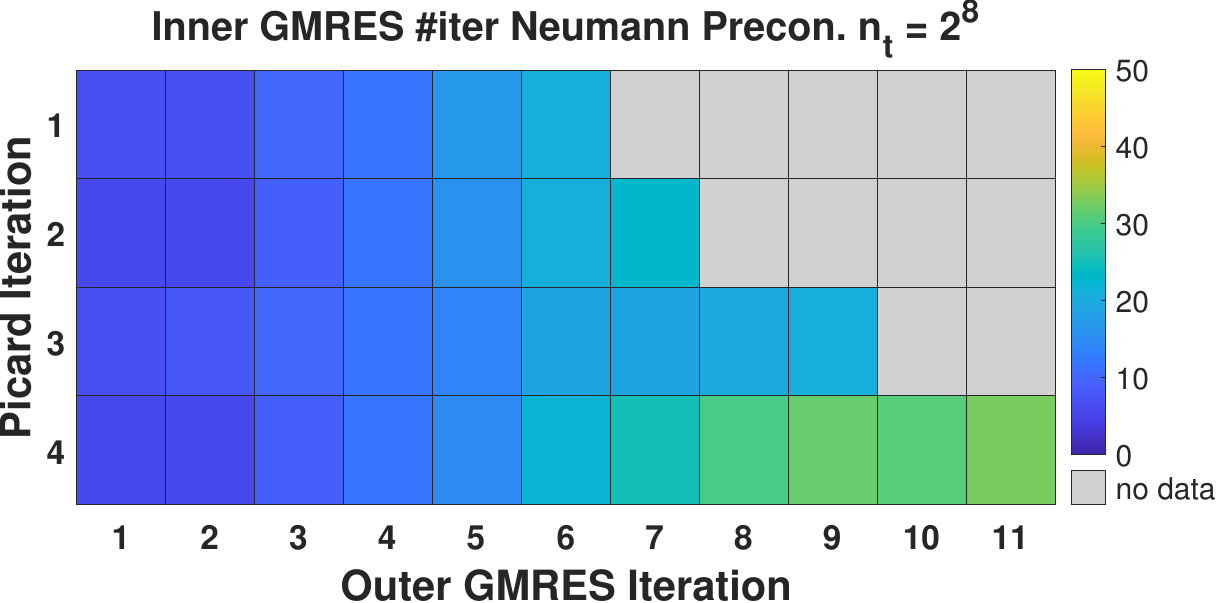}
    \hfill
\includegraphics[trim={0.8cm 0 0 0}, clip, width=0.5\textwidth]{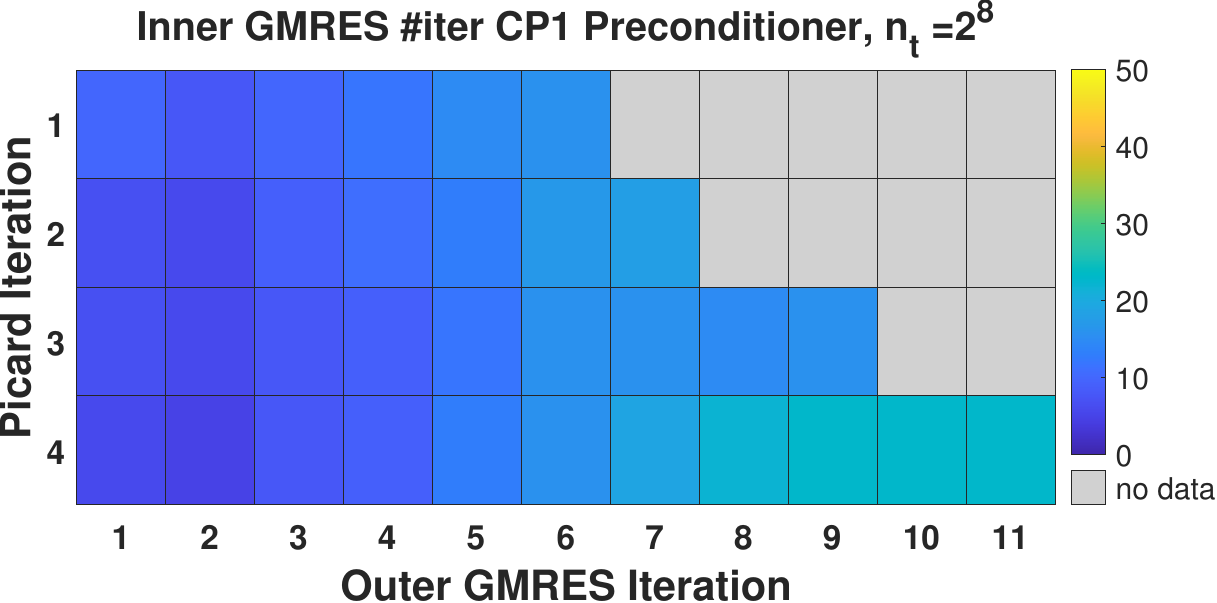}
    \caption{Heatmaps detailing the number of inner GMRES iterations nested within each outer GMRES step across the sequential Picard iterations. For uniform time-stepping with $n_t = 2^8$, and Coefficient of Variation $\sigma_\nu = 0.1$, the CP1 preconditioner (bottom right) consistently requires fewer inner iterations compared to the Mass (top left), Mean-Field (top right), and Neumann (bottom left) preconditioners.}
    \label{fig:iteration-counts-prec}
\end{figure}

\begin{figure}[ht]
    \centering
    \includegraphics[trim={1.5cm 0.8cm 2.75cm 0}, clip, width=0.47\textwidth]{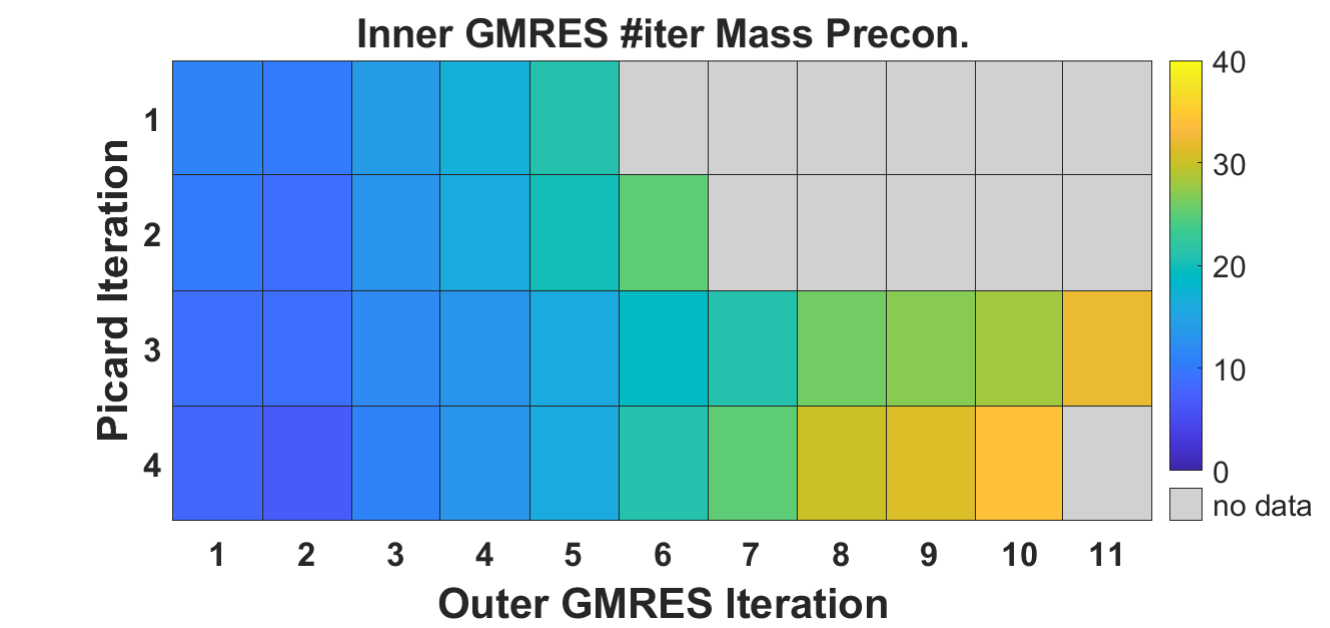}
    \hfill
    \includegraphics[trim={2.4cm 0.9cm 0 0}, clip, width=0.52\textwidth]{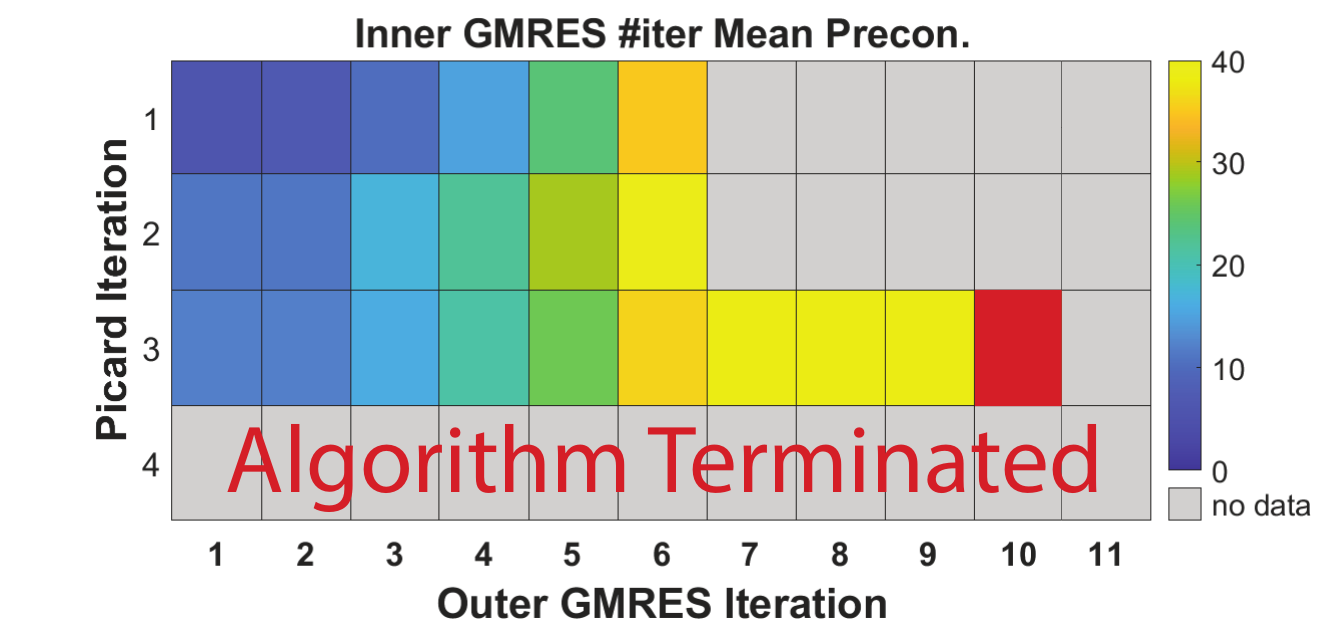}
   \\
   \vspace{0.2cm}
    \includegraphics[trim={0 0 2.5cm 0}, clip, width=0.475\textwidth]{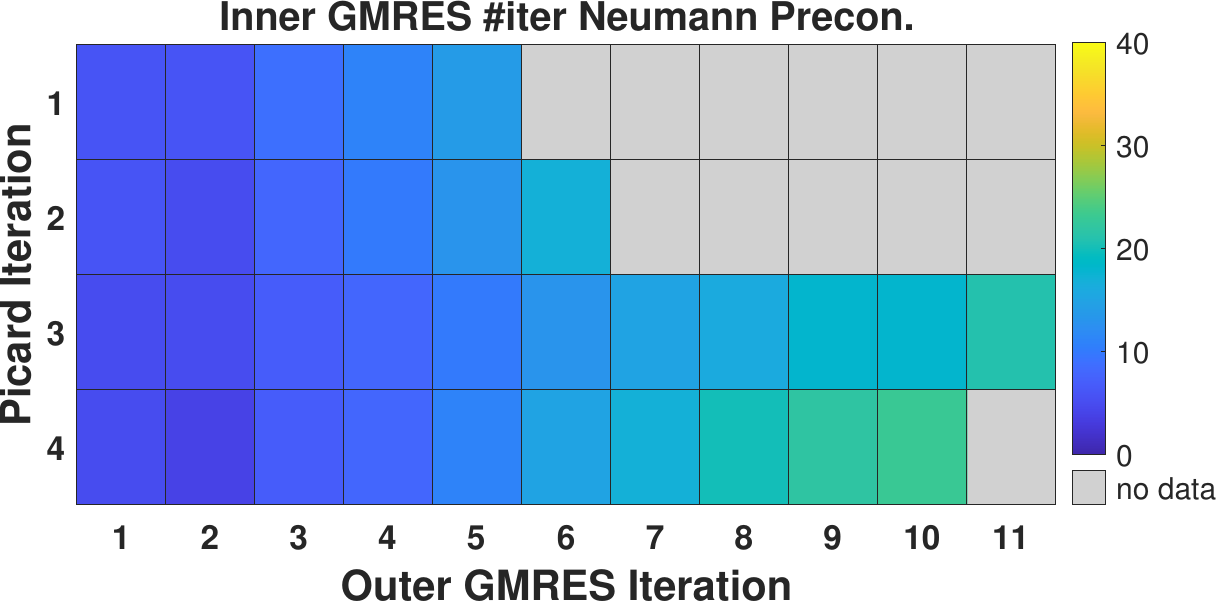}
    \hfill
\includegraphics[trim={0.8cm 0 0 0}, clip, width=0.515\textwidth]{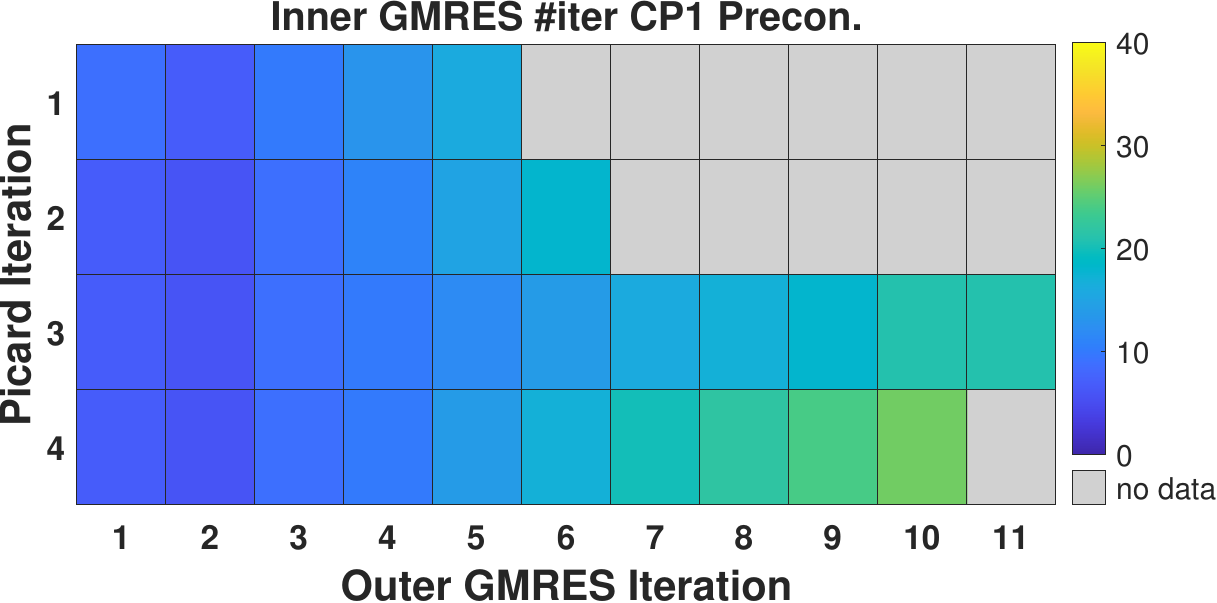}
    \caption{Heatmaps detailing the number of inner GMRES iterations nested within each outer GMRES step across the sequential Picard iterations. For nonuniform time-stepping, and Coefficient of Variation $\sigma_\nu = 0.1$, the CP1 preconditioner (bottom right) consistently requires fewer inner iterations compared to the Mass (top left), and Neumann (bottom left) preconditioners. The algorithm with the Mean-Field (top right) preconditioner terminated prematurely due to not reaching desired tolerance levels within $100$ iterations.}
    \label{fig:iteration-counts-prec-seq}
\end{figure}

\begin{figure}[H]
    \centering
    \includegraphics[trim={1cm 0cm 1cm .1cm}, clip, width=0.47\textwidth]{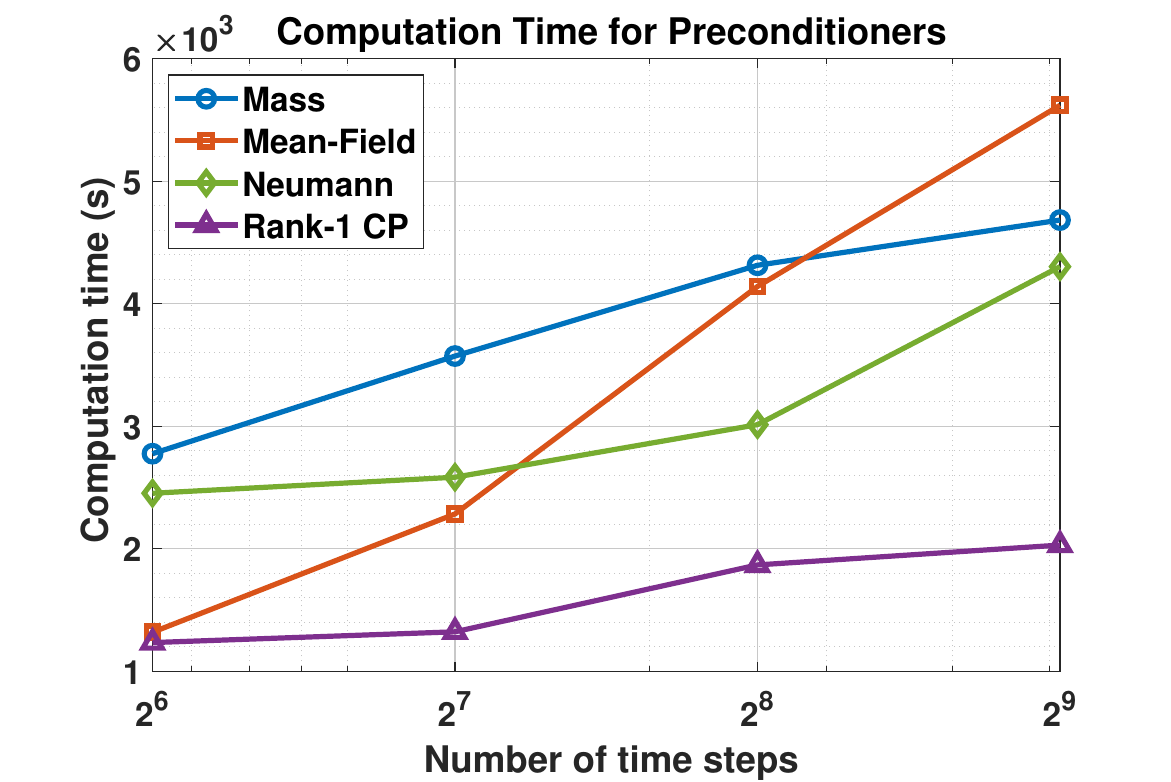}
        \includegraphics[trim={0cm 0cm 1cm .1cm}, clip, width=0.5\textwidth]{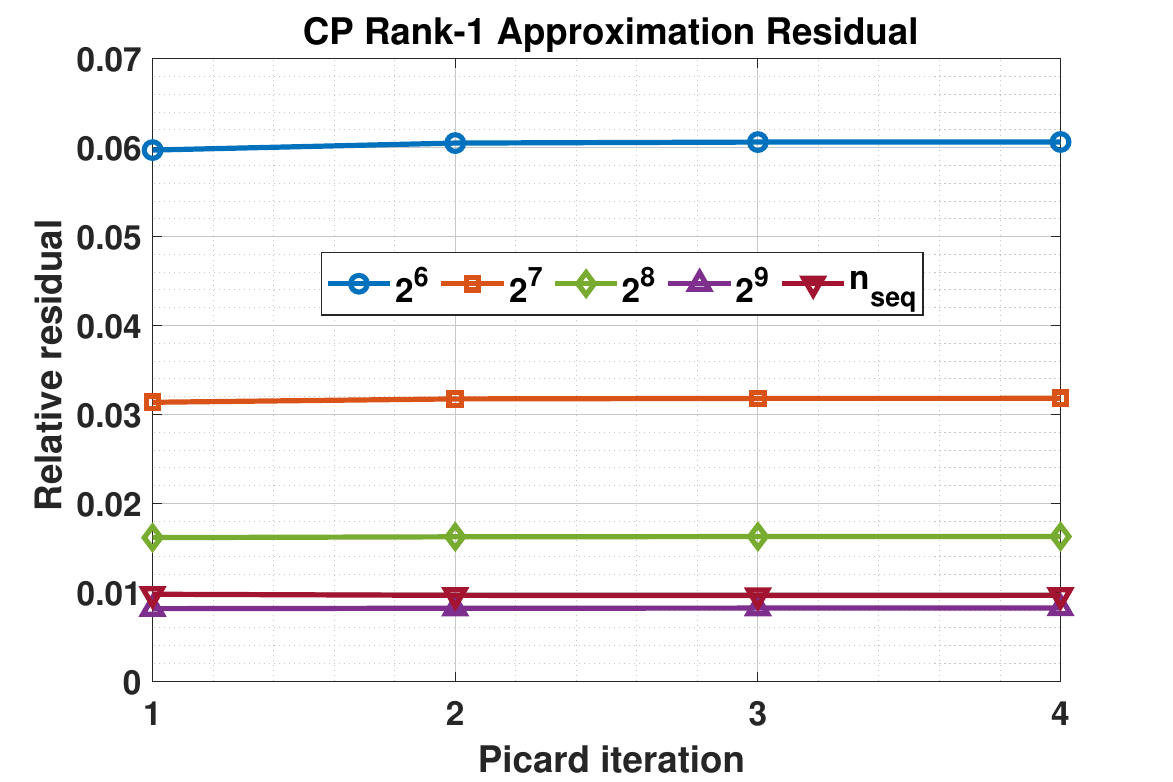}
    \caption{Left: Total computational time (CPU seconds) for each preconditioner at different time steps ($n_t$). Right: Evolution of the relative residual for the rank-one CP-ALS approximation of the global operator $\mathbb{F}_u$ across successive nonlinear Picard iterations, evaluated for varying temporal resolutions.}
    \label{fig:comptime-precon}
\end{figure}

\Cref{tab:comp_times} presents the total computational times corresponding to the different preconditioners. The nonuniform time-stepping regime offers superior resolution at the early transient phase, but it also presents a challenge for the preconditioners. Notably, the Mean-Field preconditioner \eqref{eq:prec-neumann-expansion} fails to converge within $100$ iterations and terminated early on the third Picard iteration. For uniform time-stepping, the CP1 preconditioner outperforms the others. The evaluation of computational times across different step sizes and preconditioners is further detailed in \cref{fig:comptime-precon} (left). In \cref{fig:iteration-counts-prec} shows the heatmap of iteration counts for the uniform time-stepping regime ($n_t = 2^8$) and in \cref{fig:iteration-counts-prec-seq} for the nonuniform time steps derived by a sequential solver (see \cref{fig:time-steps}). The computational cost of applying a single preconditioning step is broadly comparable across preconditioners as the main computation cost is application of the inverse of the rank-one Kronecker product matrices. The iteration counts supports the observed computational times, with the CP1 preconditioner consistently requiring fewer inner GMRES iterations across all Picard steps. Notably, in the nonuniform time-stepping, the Mean-Field approach fails to reach the tolerance level within $100$ iterations and terminated early at the tenth GMRES iteration of third Picard step.

In \cref{fig:comptime-precon} (right), the relative residuals generated by the CP-ALS algorithm (\cref{alg:ALS_CP}) when approximating the operator $\mathbb{F}_u$ for CP1 preconditioner are shown. Because a rank-one truncation is strictly enforced to keep the matrix inversion inexpensive, the relative residual cannot be reduced by increasing the rank. Note that the relative residual remains remarkably stable across successive Picard iterations. This indicates that the continuous updating of the nonlinear convection matrix $\mathbb{N}(\bm{u})$ does not degrade the quality of the rank-one surrogate. Moreover, the approximation error decreases as the number of time steps increases. As the time-step size becomes smaller, the separable temporal mass term $(\bm{T} \otimes \bm{I}_{n_\xi} \otimes \bm{M})$ becomes increasingly dominant. This makes the operator more amenable to a rank-one representation.

\section{Conclusion}
\label{sec:conclusion_1}
We developed a low-rank monolithic solver for the stochastic unsteady incompressible Navier--Stokes equations with uncertain viscosity. The stochastic Galerkin discretization leads to a large all-at-once system, which was treated using Tensor Train (TT) representations for the solution and CANDECOMP/PARAFAC (CP) approximations for selected matrices. Nonuniform time steps, informed by a sequential deterministic solve with the mean viscosity, were incorporated into the monolithic formulation. A flexible tolerance strategy was also implemented for the GMRES iterations and TT rounding to reduce unnecessary work during early nonlinear iterations.

The main contribution is the use of low-rank CP approximations in the preconditioners. In particular, the rank-one CP preconditioner provides a cheap approximation of the global velocity block and led to reduced iteration counts and computational times in the Narrow Channel benchmark. An artificially damped Neumann-series preconditioner was also considered for higher-rank CP approximations. The numerical results show that the proposed SGFEM solver agrees well with Monte Carlo and stochastic collocation reference solutions. They also show that combining TT solution representations with CP-based preconditioners can effectively reduce the cost.

\section*{Acknowledgments}
This work is dedicated to the memory of Dr. Howard C. Elman, who sadly passed away during the early stages of this project. We are deeply grateful for his invaluable insights, guidance, and foundational contributions which helped shape this research.

\bibliographystyle{plain}
\bibliography{bibliography}
\end{document}